\documentclass{amsart}

\usepackage{amsfonts,amssymb,verbatim,amsmath,amsthm,latexsym,textcomp,amscd}
\usepackage{latexsym,amsfonts,amssymb,epsfig,verbatim}
\usepackage{amsmath,amsthm,amssymb,latexsym,graphics,textcomp}
\usepackage{graphicx}
\usepackage{color}
\usepackage{url}

\input xy
\xyoption{all}

\begin{document}

\newtheorem{theorem}{Theorem}[section]
\newtheorem{prop}[theorem]{Proposition}
\newtheorem{lemma}[theorem]{Lemma}
\newtheorem{cor}[theorem]{Corollary}
\newtheorem{definition}[theorem]{Definition}
\newtheorem{conj}[theorem]{Conjecture}
\newtheorem{rmk}[theorem]{Remark}
\newtheorem{claim}[theorem]{Claim}
\newtheorem{defth}[theorem]{Definition-Theorem}

\newcommand{\boundary}{\partial}
\newcommand{\C}{{\mathbb C}}
\newcommand{\integers}{{\mathbb Z}}
\newcommand{\natls}{{\mathbb N}}
\newcommand{\ratls}{{\mathbb Q}}
\newcommand{\reals}{{\mathbb R}}
\newcommand{\proj}{{\mathbb P}}
\newcommand{\lhp}{{\mathbb L}}
\newcommand{\tube}{{\mathbb T}}
\newcommand{\cusp}{{\mathbb P}}
\newcommand\AAA{{\mathcal A}}
\newcommand\BB{{\mathcal B}}
\newcommand\CC{{\mathcal C}}
\newcommand\DD{{\mathcal D}}
\newcommand\EE{{\mathcal E}}
\newcommand\FF{{\mathcal F}}
\newcommand\GG{{\mathcal G}}
\newcommand\HH{{\mathcal H}}
\newcommand\II{{\mathcal I}}
\newcommand\JJ{{\mathcal J}}
\newcommand\KK{{\mathcal K}}
\newcommand\LL{{\mathcal L}}
\newcommand\MM{{\mathcal M}}
\newcommand\NN{{\mathcal N}}
\newcommand\OO{{\mathcal O}}
\newcommand\PP{{\mathcal P}}
\newcommand\QQ{{\mathcal Q}}
\newcommand\RR{{\mathcal R}}
\newcommand\SSS{{\mathcal S}}
\newcommand\TT{{\mathcal T}}
\newcommand\UU{{\mathcal U}}
\newcommand\VV{{\mathcal V}}
\newcommand\WW{{\mathcal W}}
\newcommand\XX{{\mathcal X}}
\newcommand\YY{{\mathcal Y}}
\newcommand\ZZ{{\mathcal Z}}
\newcommand\CH{{\CC\HH}}
\newcommand\PEY{{\PP\EE\YY}}
\newcommand\MF{{\MM\FF}}
\newcommand\RCT{{{\mathcal R}_{CT}}}
\newcommand\PMF{{\PP\kern-2pt\MM\FF}}
\newcommand\FL{{\FF\LL}}
\newcommand\PML{{\PP\kern-2pt\MM\LL}}
\newcommand\GL{{\GG\LL}}
\newcommand\Pol{{\mathcal P}}
\newcommand\half{{\textstyle{\frac12}}}
\newcommand\Half{{\frac12}}
\newcommand\Mod{\operatorname{Mod}}
\newcommand\Area{\operatorname{Area}}
\newcommand\ep{\epsilon}
\newcommand\hhat{\widehat}
\newcommand\Proj{{\mathbf P}}
\newcommand\U{{\mathbf U}}
 \newcommand\Hyp{{\mathbf H}}
\newcommand\D{{\mathbf D}}
\newcommand\Z{{\mathbb Z}}
\newcommand\R{{\mathbb R}}
\newcommand\Q{{\mathbb Q}}
\newcommand\E{{\mathbb E}}
\newcommand\til{\widetilde}
\newcommand\length{\operatorname{length}}
\newcommand\tr{\operatorname{tr}}
\newcommand\gesim{\succ}
\newcommand\lesim{\prec}
\newcommand\simle{\lesim}
\newcommand\simge{\gesim}
\newcommand{\simmult}{\asymp}
\newcommand{\simadd}{\mathrel{\overset{\text{\tiny $+$}}{\sim}}}
\newcommand{\ssm}{\setminus}
\newcommand{\diam}{\operatorname{diam}}
\newcommand{\pair}[1]{\langle #1\rangle}
\newcommand{\T}{{\mathbf T}}
\newcommand{\inj}{\operatorname{inj}}
\newcommand{\pleat}{\operatorname{\mathbf{pleat}}}
\newcommand{\short}{\operatorname{\mathbf{short}}}
\newcommand{\vertices}{\operatorname{vert}}
\newcommand{\collar}{\operatorname{\mathbf{collar}}}
\newcommand{\bcollar}{\operatorname{\overline{\mathbf{collar}}}}
\newcommand{\I}{{\mathbf I}}
\newcommand{\tprec}{\prec_t}
\newcommand{\fprec}{\prec_f}
\newcommand{\bprec}{\prec_b}
\newcommand{\pprec}{\prec_p}
\newcommand{\ppreceq}{\preceq_p}
\newcommand{\sprec}{\prec_s}
\newcommand{\cpreceq}{\preceq_c}
\newcommand{\cprec}{\prec_c}
\newcommand{\topprec}{\prec_{\rm top}}
\newcommand{\Topprec}{\prec_{\rm TOP}}
\newcommand{\fsub}{\mathrel{\scriptstyle\searrow}}
\newcommand{\bsub}{\mathrel{\scriptstyle\swarrow}}
\newcommand{\fsubd}{\mathrel{{\scriptstyle\searrow}\kern-1ex^d\kern0.5ex}}
\newcommand{\bsubd}{\mathrel{{\scriptstyle\swarrow}\kern-1.6ex^d\kern0.8ex}}
\newcommand{\fsubeq}{\mathrel{\raise-.7ex\hbox{$\overset{\searrow}{=}$}}}
\newcommand{\bsubeq}{\mathrel{\raise-.7ex\hbox{$\overset{\swarrow}{=}$}}}
\newcommand{\tw}{\operatorname{tw}}
\newcommand{\base}{\operatorname{base}}
\newcommand{\trans}{\operatorname{trans}}
\newcommand{\rest}{|_}
\newcommand{\bbar}{\overline}
\newcommand{\UML}{\operatorname{\UU\MM\LL}}
\newcommand{\EL}{\mathcal{EL}}
\newcommand{\tsum}{\sideset{}{'}\sum}
\newcommand{\tsh}[1]{\left\{\kern-.9ex\left\{#1\right\}\kern-.9ex\right\}}
\newcommand{\Tsh}[2]{\tsh{#2}_{#1}}
\newcommand{\qeq}{\mathrel{\approx}}
\newcommand{\Qeq}[1]{\mathrel{\approx_{#1}}}
\newcommand{\qle}{\lesssim}
\newcommand{\Qle}[1]{\mathrel{\lesssim_{#1}}}
\newcommand{\simp}{\operatorname{simp}}
\newcommand{\vsucc}{\operatorname{succ}}
\newcommand{\vpred}{\operatorname{pred}}
\newcommand\fhalf[1]{\overrightarrow {#1}}
\newcommand\bhalf[1]{\overleftarrow {#1}}
\newcommand\sleft{_{\text{left}}}
\newcommand\sright{_{\text{right}}}
\newcommand\sbtop{_{\text{top}}}
\newcommand\sbot{_{\text{bot}}}
\newcommand\sll{_{\mathbf l}}
\newcommand\srr{_{\mathbf r}}
\newcommand\geod{\operatorname{\mathbf g}}
\newcommand\mtorus[1]{\boundary U(#1)}
\newcommand\A{\mathbf A}
\newcommand\Aleft[1]{\A\sleft(#1)}
\newcommand\Aright[1]{\A\sright(#1)}
\newcommand\Atop[1]{\A\sbtop(#1)}
\newcommand\Abot[1]{\A\sbot(#1)}
\newcommand\boundvert{{\boundary_{||}}}
\newcommand\storus[1]{U(#1)}
\newcommand\Momega{\omega_M}
\newcommand\nomega{\omega_\nu}
\newcommand\twist{\operatorname{tw}}
\newcommand\modl{M_\nu}
\newcommand\MT{{\mathbb T}}
\newcommand\Teich{{\mathcal T}}
\renewcommand{\Re}{\operatorname{Re}}
\renewcommand{\Im}{\operatorname{Im}}

\title{On Moebius and conformal maps between boundaries of CAT(-1) spaces}

\author{Kingshook Biswas}
\address{RKM Vivekananda University, Belur Math, WB-711 202, India. Email: kingshook@rkmvu.ac.in}

\begin{abstract}
We consider Moebius and conformal homeomorphisms $f : \partial X \to \partial Y$
between boundaries of CAT(-1) spaces $X,Y$ equipped with visual
metrics. A conformal map $f$ induces a topological conjugacy of the geodesic flows
of $X$ and $Y$, which is flip-equivariant if $f$ is Moebius. We define a function
$S(f) : \partial ^2 X \to \mathbb{R}$, the {\it integrated Schwarzian} of $f$,
which measures the deviation of the topological conjugacy
from being flip-equivariant, in particular vanishing if $f$ is Moebius.
Conversely if $X,Y$ are simply connected complete manifolds
with pinched negative sectional curvatures, then $f$ is Moebius
on any open set $U \subset \partial X$ such that $S(f)$ vanishes
on $\partial^2 U$. Indeed we obtain an explicit formula for the cross-ratio distortion in terms of the integrated Schwarzian.
For such manifolds, we show that there is a
Moebius homeomorphism $f : \partial X \to \partial Y$ if and only
if there is a topological conjugacy of geodesic flows $\phi : T^1
X \to T^1 Y$ with a certain uniform continuity property along
geodesics.

\medskip

We show that if $X,Y$ are
proper, geodesically complete CAT(-1) spaces then any Moebius homeomorphism
$f$ extends to a $(1, \log 2)$-quasi-isometry with image
$\frac{1}{2}\log 2$-dense in $Y$. We prove that if $X,Y$ are in
addition metric trees then $f$ extends to a surjective
isometry. The proofs involve a study of a space
$\mathcal{M}(\partial X)$ of metrics on $\partial X$ Moebius
equivalent to a visual metric and a natural isometric embedding of $X$
into $\mathcal{M}(\partial X)$. For $C^1$ conformal maps $f : \partial X \to \partial Y$ with bounded integrated Schwarzian
and with domain $X$ a simply connected negatively curved manifold with a lower bound on sectional curvature, similar arguments
show that $f$ extends to a $(1, \log 2 + 12||S(f)||_{\infty})$ quasi-isometry.

\medskip

We also obtain a dynamical
classification of Moebius self-maps $f : \partial X \to \partial
X$ into three types, elliptic, parabolic and hyperbolic.
\end{abstract}

\bigskip

\maketitle

\tableofcontents

\section{Introduction}

\medskip

The problems we consider in this article are motivated by rigidity
results for negatively curved manifolds. The Mostow Rigidity
Theorem asserts than an isomorphism between fundamental groups of
closed hyperbolic $n$-manifolds (where $n \geq 3$) is induced by
an isometry between the manifolds. Thus hyperbolic manifolds are
determined upto isometry by their fundamental groups. It is natural to ask
for closed manifolds with variable negative curvature what extra information
over and above the fundamental group is required to determine the metric. Recall
that each free homotopy class of closed curves in a closed
negatively curved manifold contains a unique closed geodesic. Thus
a closed negatively curved manifold $X$ comes equipped with a
length function $l_X : \pi_1(X) \to \mathbb{R}^+$ (which is
constant on conjugacy classes). The marked length spectrum
rigidity problem asks whether the pair $(\pi_1(X), l_X)$ (the
{\it marked length spectrum} of $X$) determines the manifold $X$ upto isometry.
More precisely, if $X,Y$ are closed negatively curved $n$-manifolds and $\Phi : \pi_1(X) \to
\pi_1(Y)$ is an isomorphism such that $l_X = l_Y \circ \Phi$, then
is $\Phi$ induced by an isometry $F : X \to Y$?

\medskip

Otal proved that this is indeed the case if the dimension $n = 2$ \cite{otal2}.
The problem remains open in higher dimensions. It is known however
to be equivalent to two related problems, which we briefly
describe. The geodesic conjugacy problem asks whether the
existence of a homeomorphism between the unit tangent bundles
$\phi : T^1 X \to T^1 Y$ conjugating the geodesic flows
implies isometry of the manifolds. Hamenstadt proved that equality
of marked length spectra is equivalent to existence of a
geodesic conjugacy \cite{hamenstadt1}. Thus the problems of marked length spectrum
rigidity and geodesic conjugacy are equivalent.

\medskip

We recall that the boundary at infinity $\partial X$ of a CAT(-1)
space carries a natural class of metrics $\rho_x, x \in X$ called
visual metrics, which are Moebius equivalent, in the sense that
metric cross-ratios are the same for all metrics $\rho_x$. For
background on visual metrics and cross-ratios we refer to Bourdon
\cite{bourdon1}, \cite{bourdon2}. Recall that a continuous embedding $f :
\partial X \to \partial Y$ between boundaries of CAT(-1) spaces $X,Y$ is Moebius
if it preserves cross-ratios. Any isometric embedding $F : X \to
Y$ extends to a Moebius embedding $\partial F : \partial X \to
\partial Y$. Bourdon showed in \cite{bourdon1}, that for a
Gromov-hyperbolic group $\Gamma$ with two quasi-convex actions on
CAT(-1) spaces $X, Y$, the natural $\Gamma$-equivariant
homeomorphism $f$ between the limit sets $\Lambda X, \Lambda Y$ is Moebius if and only if
there is a $\Gamma$-equivariant conjugacy of the abstract geodesic
flows $\mathcal{G}\Lambda X$ and $\mathcal{G} \Lambda Y$ compatible with
$f$. In particular for $X, Y$ the universal covers of two closed
negatively curved manifolds $\overline{X},\overline{Y}$, it follows that
the geodesic flows of $\overline{X},\overline{Y}$ are conjugate if
and only if the induced equivariant boundary map $f : \partial X
\to \partial Y$ is Moebius.

\medskip

Otal showed that equality of the marked length spectrum for two negatively curved metrics
on the same closed manifold is equivalent to the existence of an equivariant Moebius map between the
boundaries at infinity of the universal covers \cite{otal1}. We remark that the same conclusion holds when
the marked length spectra of two closed negatively curved
manifolds coincide (the manifolds are not necessarily assumed to
be diffeomorphic), using the results of Hamenstadt (equality of the marked length spectrum being
equivalent to conjugacy of geodesic flows) and Bourdon (conjugacy of geodesic flows being equivalent to the
boundary map being Moebius).

\medskip

It follows that the marked length spectrum, the geodesic flow, and the
Moebius structure on the boundary at infinity of the universal cover are all equivalent
data for a closed negatively curved manifold, and the question becomes whether any one of these is enough
to determine the metric. We discuss in section 5 the proofs of
these equivalences. In the case of
simply connected, complete Riemannian manifolds of sectional
curvature bounded above by $-1$, the marked length spectrum no longer makes sense, but
one may still consider the correspondence between Moebius maps and geodesic
conjugacies. We define a certain uniform continuity property for geodesic conjugacies,
{\it uniform continuity along geodesics} (which is satisfied in particular by uniformly continuous
maps). Recall that a CAT(-1) space $X$ is {\it
geodesically complete} if every geodesic segment in $X$ can be
extended (not necessarily uniquely) to a bi-infinite geodesic. We
show in section 4:

\medskip

\begin{theorem} \label{geodesicmoebius} Let $X$ be a simply
connected complete Riemannian manifold with sectional curvatures
bounded above by $-1$, and let $Y$ be a proper geodesically complete CAT(-1) space.
If there is a homeomorphism $\phi : T^1 X \to \mathcal{G}Y$
conjugating the geodesic flows of $X$ and $Y$ which is uniformly continuous along geodesics
then $\phi$ induces a map $f : \partial X \to \partial Y$ which is Moebius.
\end{theorem}

\medskip

Recall that there is a notion of a conformal homeomorphism between
metric spaces, in particular between boundaries of CAT(-1) spaces equipped
with visual metrics. We consider {\it $C^1$ conformal} maps, i.e. those for which
the pointwise derivative is a continuous function. A $C^1$ conformal map
$f : \partial X \to \partial Y$ between boundaries of CAT(-1) spaces
induces a topological conjugacy $\phi : \mathcal{G}X \to \mathcal{G}Y$
between the abstract geodesic flows
of $X$ and $Y$ (following Bourdon \cite{bourdon1}), where $\mathcal{G}X, \mathcal{G}Y$
are the spaces of bi-infinite geodesics in $X$ and $Y$. The conjugacy is
equivariant with respect to the flips if $f$ is Moebius.
We define a function $S(f) : \partial^2 X \to \mathbb{R}$, the {\it integrated Schwarzian} of
$f$, which measures the deviation of the conjugacy from being flip-equivariant, vanishing
in particular if $f$ is Moebius. Coversely, if the domain $X$ is a simply connected negatively curved manifold
also satisfying a lower curvature bound $-b^2 \leq K \leq -1$, then, as in the classical case, bounds on the integrated Schwarzian imply bounds
on cross-ratio distortion. Indeed we have an exact formula for the cross-ratio distortion:

\medskip

\begin{theorem} \label{crossratiodistn} Let $X$ be a simply
connected complete Riemannian manifold with sectional curvatures
satisfying $-b^2 \leq K \leq -1$ for some $b \geq 1$, and let $Y$ be a proper geodesically complete CAT(-1) space.
Let $f : U \subset \partial X \to V \subset \partial Y$ be a $C^1$ conformal map between open subsets
$U, V$. Then
$$
\log \frac{[f(\xi), f(\xi'), f(\eta), f(\eta')]}{[\xi, \xi', \eta, \eta']} = \frac{1}{2}\left(S(f)(\xi,\eta) + S(f)(\xi',\eta') - S(f)(\xi, \eta') - S(f)(\xi', \eta)\right)
$$
for all $(\xi, \xi',\eta,\eta') \in \partial^4 U$.
\end{theorem}

\medskip

The integrated Schwarzian also satisfies a cocycle identity, thus two $C^1$ conformal maps differ by post-composition with a
Moebius map if and only if their integrated Schwarzians are equal, as in the classical case.

\medskip

\medskip

In the case of a lower curvature bound we
have a converse to Theorem \ref{geodesicmoebius} above:

\medskip

\begin{theorem} \label{mobconfgeo} Let $X,Y$ be as in the previous theorem and $f : \partial X \to \partial Y$ a $C^1$ conformal map.
Then the induced topological conjugacy of geodesic flows $\phi : T^1 X \to \mathcal{G}Y$ is uniformly continuous along
geodesics if and only if $f$ is Moebius.
\end{theorem}

\medskip

We then consider in the more general context of CAT(-1) spaces, the question of whether a Moebius
embedding $f : \partial X \to \partial Y$ between the boundaries
of two CAT(-1) spaces extends to an isometric embedding $F : X \to
Y$. In \cite{bourdon2}, Bourdon proved the following Theorem:

\medskip

\begin{theorem} \label{bourdonthm} {\bf (Bourdon)} If $X$ is a rank one symmetric
space and $Y$ a CAT(-1) space then any Moebius embedding $f :
\partial X \to \partial Y$ extends to an isometric embedding $F :
X \to Y$.
\end{theorem}

\medskip

We consider the general case where the domain $X$ is an
arbitrary CAT(-1) space. We
prove the following in section 6:

\medskip

\begin{theorem} \label{mainthm1} Let $X, Y$ be proper geodesically complete CAT(-1)
spaces such that $\partial X$ has at least four points, and let
$f : \partial X \to \partial Y$ be a Moebius homeomorphism. Then $f$ extends to a $(1, \log
2)$-quasi-isometry $F : X \to Y$, with image $\log 2$-dense in $Y$.
\end{theorem}

\medskip

In the case of metric trees we have:

\medskip

\begin{theorem} \label{mainthm2} Let $X, Y$ be proper geodesically
complete metric trees such that $\partial X$ has at least four points and
let $f : \partial X \to \partial Y$ be a Moebius homeomorphism.
Then $f$ extends to a surjective isometry $F : X \to Y$.
\end{theorem}

\medskip

The proofs of the above Theorems involve a study of the space
$\mathcal{M}(\partial X)$ of metrics on the boundary $\partial X$ of a
proper geodesically complete CAT(-1) space $X$ which are Moebius
equivalent to a visual metric. The key point is that
there is a natural metric $d_{\mathcal{M}}$ on $\mathcal{M}(\partial X)$ such that
the map $i_X : X \to \mathcal{M}(\partial X)$ sending a point $x \in X$
to the visual metric $\rho_x$ based at $x$ is an isometric
embedding. The space $(\mathcal{M}(\partial X), d_{\mathcal{M}})$ is itself isometric to a closed, locally compact
subspace of the Banach space $C(\partial X)$ of continuous functions on $\partial X$.
By studying the derivative of the embedding $i_X$ along geodesics in $X$, we show that it has image
$\frac{1}{2}\log 2$-dense in $\mathcal{M}(\partial X)$, and is surjective in the case of a metric tree. Thus
we may define a nearest-point projection map (not unique) $\pi_X : \mathcal{M}(\partial X) \to X$ which is a $(1,\log 2)$
quasi-isometry.

\medskip

A Moebius map $f : \partial X \to \partial Y$ induces a natural map $\hat{f} : \mathcal{M}(\partial X) \to \mathcal{M}(\partial Y)$
(by push-forward of metrics) which is a surjective isometry. The extension $F : X \to Y$ of $f$ is then defined by $F = \pi_Y \circ \hat{f} \circ i_X$.

\medskip

For $C^1$ conformal maps with bounded integrated Schwarzian, similar arguments lead to the following:

\medskip

\begin{theorem}\label{confextn} Let $X$ be a simply
connected complete Riemannian manifold with sectional curvatures
satisfying $-b^2 \leq K \leq -1$ for some $b \geq 1$, and let $Y$ be a proper geodesically complete CAT(-1) space.
Let $f : \partial X \to \partial Y$ be a $C^1$ conformal map such that $S(f)$ is bounded. Then $f$ extends to a $(1, \log
2 + 12||S(f)||_{\infty})$-quasi-isometry $F : X \to Y$. If $Y$ is also a simply
connected complete Riemannian manifold with sectional curvatures
satisfying $-b^2 \leq K \leq -1$ for some $b \geq 1$, then the image is ($\log 2 + 12||S(f)||_{\infty}$)-dense in $Y$.
\end{theorem}

We have as corollaries of the above theorems the
following:

\medskip

\begin{theorem} \label{geodconj} Let $X$ be a simply
connected complete Riemannian manifold with sectional curvatures
bounded above by $-1$ and $Y$ a proper geodesically complete CAT(-1) space. Suppose that there is a conjugacy
$\phi : T^1 X \to \mathcal{G}Y$ of geodesic flows which is uniformly continuous along geodesics. Then:

\smallskip

{\noindent (1)} There is a $(1,\log 2)$-quasi-isometry $F : X \to Y$ with image
$\frac{1}{2}\log 2$-dense in $Y$.

\smallskip

{\noindent (2)} If $X$ is a rank one symmetric space then $F$ can
be taken to be a surjective isometry.
\end{theorem}

\medskip

We remark that part (2) of the above theorem implies as a corollary marked length
spectrum rigidity for rank one locally symmetric spaces, a
well-known fact proved earlier by Hamenstadt \cite{hamenstadt2} using the
celebrated minimal entropy rigidity theorem of Besson-Courtois-Gallot \cite{bcg1}.

\medskip

Finally we obtain in section 7 a dynamical classification of Moebius self-maps
into three types, elliptic, parabolic and hyperbolic:

\medskip

\begin{theorem} \label{classfn} Let $X$ be a proper geodesically
complete CAT(-1) space and $f : \partial X \to \partial X$ a
Moebius self-map of its boundary. Then one of the following three mutually
exclusive cases holds:

\smallskip

\noindent (1) For all $x \in X$, the iterates $f^n : (\partial X,
\rho_x) \to (\partial X, \rho_x)$ are uniformly bi-Lipschitz (we
say $f$ is {\it elliptic}).

\smallskip

\noindent (2) There is a unique fixed point $\xi_0 \in \partial X$
of $f$ such that $f^n(\xi) \to \xi_0$ for all $\xi$ as $n \to \pm
\infty$ (we say $f$ is {\it parabolic}).

\smallskip

\noindent (3) There is a pair of distinct fixed points $\xi_+,
\xi_-$ of $f$ such that for all $\xi \in \partial X - \{ \xi_+,
\xi_- \}$, $f^n(\xi) \to \xi_+$ as $n \to +\infty$ and $f^n(\xi)
\to \xi_-$ as $n \to -\infty$ (we say $f$ is {\it hyperbolic}).
\end{theorem}

\medskip

\noindent {\bf Acknowledgements.} The author is grateful to Marc
Bourdon and Mahan Mj for helpful discussions. The author was
supported by CEFIPRA grant no. 4301-1: "Kleinian groups: geometric
and analytic aspects".

\medskip

\section{Spaces of Moebius equivalent metrics}

\medskip

Let $(Z,\rho_0)$ be a compact metric space with at least four points. For a metric $\rho$ on
$Z$ we define the metric cross-ratio with respect to $\rho$ of a quadruple of distinct
points $(\xi, \xi', \eta, \eta')$ of $Z$ by
$$
[\xi \xi' \eta \eta']_{\rho} := \frac{\rho(\xi, \eta) \rho(\xi', \eta')}{\rho(\xi,
\eta')\rho(\xi', \eta)}
$$
We say that a diameter one metric $\rho$ on $Z$ is {\it antipodal} if for
any $\xi \in Z$ there exists $\eta \in Z$ such that $\rho(\xi,
\eta) = 1$. We assume that $\rho_0$ is diameter one and antipodal. We say two metrics $\rho_1, \rho_2$
on $Z$ are {\it Moebius equivalent} if their metric cross-ratios agree:
$$
[\xi \xi' \eta \eta']_{\rho_1} = [\xi \xi' \eta \eta']_{\rho_2}
$$
for all $(\xi, \xi', \eta, \eta')$. We define
$$
\mathcal{M}(Z, \rho_0) := \{ \rho : \rho \hbox{ is an antipodal, diameter one metric on } Z
\hbox{ Moebius equivalent to } \rho_0 \}
$$
We will write $\mathcal{M}(Z, \rho_0) = \mathcal{M}$. Note we do not assume
that the metrics $\rho \in \mathcal{M}$ induce the same topology on $Z$ as $\rho_0$,
but we will show that they are indeed all bi-Lipschitz equivalent to each other.
For $\rho_1, \rho_2 \in \mathcal{M}$
we define a
positive function on $Z$ called the derivative of $\rho_2$ with respect to
$\rho_1$ by
$$
\frac{d\rho_2}{d\rho_1}(\xi) := \frac{\rho_2(\xi,
\eta)\rho_2(\xi,\eta')\rho_1(\eta, \eta')}{\rho_1(\xi,
\eta)\rho_1(\xi,\eta')\rho_2(\eta, \eta')}
$$
where $\eta, \eta' \in Z$ are distinct points not equal to $\xi$.

\medskip

\begin{lemma} \label{welldefined} The function
$\frac{d\rho_2}{d\rho_1}$ is well-defined.
\end{lemma}

\medskip

\noindent{\bf Proof:} Given two pairs of distinct points $\eta,
\eta'$ and $\beta, \beta'$ not equal to $x$, the desired equality
$$
\frac{\rho_2(\xi,
\eta)\rho_2(\xi,\eta')\rho_1(\eta, \eta')}{\rho_1(\xi,
\eta)\rho_1(\xi,\eta')\rho_2(\eta, \eta')} = \frac{\rho_2(\xi,
\beta)\rho_2(\xi,\beta')\rho_1(\beta, \beta')}{\rho_1(\xi,
\beta)\rho_1(\xi,\beta')\rho_2(\beta, \beta')}
$$
follows from the equality
$$
[\xi \beta \eta \eta']_{\rho_2} [\xi \eta \eta' \beta']_{\rho_2} =
[\xi \beta \eta \eta']_{\rho_1} [\xi \eta \eta' \beta']_{\rho_1}
$$
$\diamond$

\medskip

The next Lemma follows from a straightforward computation using
the definition of the derivative, we omit the proof:

\medskip

\begin{lemma} \label{chainrule} {\bf (Chain Rule)} For $\rho_1,
\rho_2, \rho_3 \in \mathcal{M}$ we have
$$
\frac{d\rho_3}{d\rho_1} = \frac{d\rho_3}{d\rho_2} \frac{d\rho_2}{d\rho_1}
$$
and
$$
\frac{d\rho_2}{d\rho_1} = 1/\left(\frac{d\rho_1}{d\rho_2}\right)
$$
\end{lemma}

\medskip

\begin{lemma} \label{bounded} For $\rho \in \mathcal{M}$ the function
$f = \frac{d\rho}{d\rho_0}$ is bounded.
\end{lemma}

\medskip

\noindent{\bf Proof:} Suppose not, let $\xi_n \in Z$ be a sequence such that
$f(\xi_n) \to \infty$. Passing to a subsequence we may assume $\xi_n \to \xi$,
choose $\eta, \eta'$ distinct points in $Z$ not equal to $\xi$, then we have
\begin{align*}
\limsup f(\xi_n) & = \limsup \frac{\rho(\xi_n,
\eta)\rho(\xi_n,\eta')\rho_0(\eta,\eta')}{\rho_0(\xi_n,
\eta)\rho_0(\xi_n,\eta')\rho(\eta,\eta')} \\
& \leq \frac{1}{\rho_0(\xi,
\eta)\rho_0(\xi,\eta')\rho(\eta,\eta')}, \\
\end{align*}
a contradiction. $\diamond$

\medskip

\begin{lemma} \label{mvt} {\bf (Geometric Mean-Value Theorem)}
$$
\rho_2(\xi, \eta)^2 = \rho_1(\xi, \eta)^2
\frac{d\rho_2}{d\rho_1}(\xi) \frac{d\rho_2}{d\rho_1}(\eta)
$$
for all $\xi, \eta \in Z$.
\end{lemma}

\medskip

\noindent{\bf Proof:} Given $\xi \neq \eta$ choose a point $\beta$
distinct from $\xi, \eta$, then by definition we may write
$$
\frac{d\rho_2}{d\rho_1}(\xi) =
\frac{\rho_2(\xi,\eta)\rho_2(\xi,\beta)\rho_1(\eta,\beta)}{\rho_1(\xi,\eta)\rho_1(\xi,\beta)\rho_2(\eta,\beta)}
\ , \ \frac{d\rho_2}{d\rho_1}(\eta) =
\frac{\rho_2(\eta,\xi)\rho_2(\eta,\beta)\rho_1(\xi,\beta)}{\rho_1(\eta,\xi)\rho_1(\eta,\beta)\rho_2(\xi,\beta)}
$$
from which it follows that
$$
\frac{d\rho_2}{d\rho_1}(\xi)\frac{d\rho_2}{d\rho_1}(\eta) =
\left(\frac{\rho_2(\xi, \eta)}{\rho_1(\xi, \eta)}\right)^2
$$
$\diamond$

\medskip

For $\rho \in \mathcal{M}$ since $\frac{d\rho}{d\rho_0}$ is bounded it
follows from the above Lemma that $\rho \leq K \rho_0$, hence the
functions $\xi \mapsto \rho(\xi,\eta)$ are continuous for
all $\eta \in Z$, therefore the functions $\frac{d\rho}{d\rho_0}$ are
continuous. Since $\frac{d\rho_2}{d\rho_1} = \frac{d\rho_2}{d\rho_0}/\frac{d\rho_1}{d\rho_0}$
it follows that all functions $\frac{d\rho_2}{d\rho_1}$ are continuous, so bounded above and below
by positive constants, hence by the above Lemma
all metrics $\rho \in \mathcal{M}$ are bi-Lipschitz to each other
and induce the same topology on $Z$ as $\rho_0$.
The following Lemma justifies the
use of the term 'derivative':

\medskip

\begin{lemma} \label{derivative} If $\xi \in Z$ is not an isolated
point then
$$
\frac{d\rho_2}{d\rho_1} = \lim_{\eta \to \xi} \frac{\rho_2(\xi,
\eta)}{\rho_1(\xi, \eta)}
$$
\end{lemma}

\medskip

\noindent{\bf Proof:} We have
\begin{align*}
\frac{\rho_2(\xi, \eta)}{\rho_1(\xi,\eta)} & = \frac{d\rho_2}{d\rho_1}(\xi)^{1/2} \frac{d\rho_2}{d\rho_1}(\eta)^{1/2} \\
                                           & \to \frac{d\rho_2}{d\rho_1}(\xi)
                                           \\
\end{align*}
as $\eta \to \xi$. $\diamond$

\medskip

\begin{lemma} \label{maxmin}
$$
\max_{\xi \in Z} \frac{d\rho_2}{d\rho_1}(\xi) \cdot \min_{\xi \in Z}
\frac{d\rho_2}{d\rho_1}(\xi) = 1
$$
\end{lemma}

\medskip

\noindent{\bf Proof:} Let $\lambda, \mu$ denote the maximum and
minimum values of $\frac{d\rho_2}{d\rho_1}$ respectively, and let $\xi, \eta \in Z$ denote
points where the maximum and minimum values are attained
respectively. Choosing $\eta' \in Z$ such that $\rho_1(\xi, \eta') =
1$ gives
$$
1 \geq \rho_2(\xi, \eta') = \frac{d\rho_2}{d\rho_1}(\xi)^{1/2}
\frac{d\rho_2}{d\rho_1}(\eta')^{1/2} \geq \lambda^{1/2} \cdot
\mu^{1/2}
$$
while choosing $\xi' \in Z$ such that $\rho_2(\xi', \eta) = 1$
gives
$$
1 \geq \rho_1(\xi', \eta) =
1/\left(\frac{d\rho_2}{d\rho_1}(\xi')^{1/2}
\frac{d\rho_2}{d\rho_1}(\eta)^{1/2}\right) \geq 1/(\lambda^{1/2}
\mu^{1/2})
$$
hence $\lambda \cdot \mu = 1$.$\diamond$

\medskip

We now define for $\rho_1, \rho_2 \in \mathcal{M}$,
$$
d_{\mathcal{M}}(\rho_1, \rho_2) := \max_{\xi \in Z}
\log \frac{d\rho_2}{d\rho_1}(\xi)
$$

\medskip

\begin{lemma} \label{metric} The function $d_{\mathcal{M}}$
is a metric on $\mathcal{M}$.
\end{lemma}

\medskip

\noindent{\bf Proof:} For $\rho_1, \rho_2 \in \mathcal{M}$, $(\max_{\xi \in Z}
\frac{d\rho_2}{d\rho_1}(\xi))^2 \geq (\max_{\xi \in Z}
\frac{d\rho_2}{d\rho_1}(\xi)) \cdot (\min_{\xi \in Z}
\frac{d\rho_2}{d\rho_1}(\xi)) = 1$, hence $d_{\mathcal{M}}(\rho_1,
\rho_2) \geq 0$. Moreover $d_{\mathcal{M}}(\rho_1,
\rho_2) = 0$ implies $\max_{\xi \in Z}
\frac{d\rho_2}{d\rho_1}(\xi) = 1$ hence $\min_{\xi \in Z}
\frac{d\rho_2}{d\rho_1}(\xi) = 1$ by the previous Lemma, hence
$\frac{d\rho_2}{d\rho_1} \equiv 1$, and it then follows from the
Geometric Mean-Value Theorem that $\rho_1 \equiv \rho_2$.

\medskip

Symmetry of $d_{\mathcal{M}}$ follows from
$\frac{d\rho_1}{d\rho_2} = 1/\frac{d\rho_2}{d\rho_1}$ and the
previous Lemma, while the triangle inequality follows easily from
the Chain Rule $\frac{d\rho_3}{d\rho_1} = \frac{d\rho_3}{d\rho_2}
\frac{d\rho_2}{d\rho_1}$. $\diamond$

\medskip

Let $(C(Z), ||\cdot||_{\infty})$ denote the Banach space of
continuous functions on $Z$ equipped with the supremum norm.

\medskip

\begin{lemma} \label{cz} The map

\begin{align*}
\mathcal{M} & \to C(Z) \\
   \rho     & \mapsto \log \frac{d\rho}{d\rho_0} \\
\end{align*}

is an isometric embedding.
\end{lemma}

\medskip

\noindent{\bf Proof:} It follows from Lemma \ref{maxmin} that
$\max_{\xi \in Z} \log \frac{d\rho_2}{d\rho_1}(\xi) = ||\log \frac{d\rho_2}{d\rho_1}||_{\infty}$,
hence
$$
d_{\mathcal{M}}(\rho_1, \rho_2) = \left|\left|\log
\frac{d\rho_2}{d\rho_1}\right|\right|_{\infty} = \left|\left|\log
\frac{d\rho_2}{d\rho_0} - \log \frac{d\rho_1}{d\rho_0}\right|\right|_{\infty}
$$
(where the second equality uses the Chain Rule). $\diamond$

\medskip

\begin{lemma} \label{closed} The image of the above embedding is
closed in $C(Z)$.
\end{lemma}

\medskip

\noindent{\bf Proof:} Let $\rho_n \in \mathcal{M}$ such that $g_n
= \log \frac{d\rho_n}{d\rho_0}$ converges in $C(Z)$ to $g$. Define
$f = e^g$ and $\rho(\xi, \eta) := \rho_0(\xi,
\eta)f(\xi)^{1/2}f(\eta)^{1/2}, \xi, \eta \in Z$, then it follows from the
Geometric Mean Value Theorem that $\rho(\xi, \eta) = \lim \rho_n(\xi, \eta)$. Passing to the limit
in the triangle inequality for $\rho_n$ gives the triangle inequality for
$\rho$, while symmetry and positivity of $\rho$ are clear, hence $\rho$
is a metric. Moreover it follows easily from the definition of
$\rho$ that $\rho$ is Moebius equivalent to $\rho_0$, and moreover $\frac{d\rho}{d\rho_0} = f$.
Since the $\rho_n$'s have diameter one it follows that $\rho$ has diameter
less than or equal to one. Given $\xi \in Z$ let $\eta_n \in Z$
such that $\rho_n(\xi, \eta_n) = 1$, passing to a subsequence we may assume
$\eta_n$ converges to some $\eta$, then
\begin{align*}
|\rho(\xi, \eta) - \rho_n(\xi,\eta_n)| & \leq |\rho(\xi, \eta) -
\rho_n(\xi, \eta)| + |\rho_n(\xi,\eta) - \rho_n(\xi, \eta_n)| \\
& \leq |\rho(\xi, \eta) - \rho_n(\xi, \eta)| + \rho_n(\eta,
\eta_n) \\
 & \to 0 \\
\end{align*}
since $\rho_0(\eta, \eta_n) \to 0$ and the $\rho_n$'s are
uniformly bi-Lipschitz equivalent to $\rho_0$ (being a bounded
sequence in $\mathcal{M}$), hence $\rho(\xi, \eta) = 1$. Thus
$\rho$ is of diameter one and is antipodal, hence $\rho \in
\mathcal{M}$ and $g$ is the image of $\rho$ under the isometric
embedding. $\diamond$

\medskip

\begin{lemma} \label{lipschitz} The function $f = \frac{d\rho_2}{d\rho_1}
: (Z, \rho_1) \to \mathbb{R}$ is $K$-Lipschitz where $K = 2 (\max_{\xi \in Z} f(\xi))^2$.
\end{lemma}

\medskip

\noindent{\bf Proof:} Let $\lambda = \max_{\xi \in Z} f(\xi)$.
Let $\xi_1, \xi_2 \in Z$. We may assume $f(\xi_1) \geq f(\xi_2)$.
Choose $\xi \in Z$ such that $\rho_1(\xi_1, \xi) = 1$, then
the inequality $|\rho_2(\xi,\xi_1) - \rho(\xi,
\xi_2)| \leq \rho_2(\xi_1, \xi_2)$ gives, using the Geometric
Mean-Value Theorem,
$$
f(\xi)^{1/2} \left|f(\xi_1)^{1/2} - \rho_1(\xi,
\xi_2)f(\xi_2)^{1/2}\right| \leq \rho_1(\xi_1, \xi_2)f(\xi_1)^{1/2}f(\xi_2)^{1/2}
$$
and we have
$$
|f(\xi_1)^{1/2} - \rho_1(\xi, \xi_2) f(\xi_2)^{1/2}| = f(\xi_1)^{1/2} - \rho_1(\xi, \xi_2)
f(\xi_2)^{1/2} \geq f(\xi_1)^{1/2} - f(\xi_2)^{1/2}
$$
which, combined with the previous inequality, gives
$$
(1/\lambda^{1/2}) (f(\xi_1)^{1/2} - f(\xi_2)^{1/2}) \leq \rho_1(\xi_1,
\xi_2) \lambda
$$
hence
\begin{align*}
|f(\xi_1) - f(\xi_2)|  & = |(f(\xi_1)^{1/2} - f(\xi_2)^{1/2})(f(\xi_1)^{1/2} +
f(\xi_2)^{1/2})| \\ & \leq \lambda^{3/2} \rho_1(\xi_1, \xi_2) 2
\lambda^{1/2} = 2 \lambda^2 \rho_1(\xi_1, \xi_2) \\
\end{align*}
$\diamond$

\medskip

\begin{lemma} \label{proper} The space $(\mathcal{M},
d_{\mathcal{M}})$ is proper, i.e. closed balls are compact. Hence $(\mathcal{M},
d_{\mathcal{M}})$ is also complete.
\end{lemma}

\medskip

\noindent{\bf Proof:} It follows from the previous Lemma that for a sequence $\rho_n \in
\mathcal{M}$ with $d_{\mathcal{M}}(\rho_n, \rho_0)$ bounded, the
functions $f_n = \frac{d\rho_n}{d\rho_0}$ are uniformly Lipschitz,
and uniformly bounded away from $0$ and $\infty$, hence the
functions $g_n = \log f_n$ are uniformly Lipschitz and uniformly bounded.
Therefore $g_n$ has a subsequence $g_{n_k}$ converging uniformly to
a continuous function $g$, which by Lemma \ref{closed} is equal to $\log
\frac{d\rho}{d\rho_0}$ for some $\rho \in \mathcal{M}$. It follows
from Lemma \ref{cz} that $\rho_{n_k} \to \rho$ in $\mathcal{M}$.
$\diamond$

\bigskip

\section{Visual metrics on the boundary of a CAT(-1) space}

\medskip

Let $(X, d_X)$ be a proper CAT(-1) space such that $\partial X$ has at least four points.

\medskip

\subsection{Definitions}

\medskip

We recall below the definitions and some
elementary properties of visual metrics and Busemann functions;
for proofs we refer to \cite{bourdon1}:

\medskip

Let $x \in X$ be a basepoint. The {\it Gromov product} of two
points $\xi, \xi' \in \partial X$ with respect to $x$ is defined by
$$
(\xi | \xi')_x = \lim_{(a,a') \to (\xi, \xi')}
\frac{1}{2}(d_X(x,a) + d_X(x,a') - d_X(a,a'))
$$
where $a,a'$ are points of $X$ which converge radially towards
$\xi$ and $\xi'$ respectively. The {\it visual metric} on
$\partial X$ based at the point $x$ is defined by
$$
\rho_x(\xi, \xi') := e^{-(\xi|\xi')_x}
$$
The distance $\rho_x(\xi,\xi')$ is less than or equal to one, with
equality iff $x$ belongs to the geodesic $(\xi \xi')$.

\medskip

\begin{lemma} \label{visualantipodal} If $X$ is geodesically
complete then $\rho_x$ is a diameter one antipodal metric.
\end{lemma}

\medskip

\noindent{\bf Proof:} Let $\xi \in \partial X$, then the geodesic
ray $[x, \xi)$ extends to a bi-infinite geodesic $(\xi' \xi)$ for
some $\xi' \in \partial X$, hence $\rho_x(\xi, \xi') = 1$, hence
$\rho_x$ is diameter one and antipodal. $\diamond$

\medskip

The Busemann function $B : \partial X \times X \times X \to
\mathbb{R}$ is defined by
$$
B(\xi, x, y) := \lim_{a \to \xi} d_X(x,a) - d_X(y,a)
$$
where $a \in X$ converges radially towards $\xi$.

\medskip

It will be convenient to consider the functions on $\partial X$, $f_{x,y}(\xi) :=
e^{B(\xi, x,y)}, g_{x,y}(\xi) = B(\xi, x, y), \xi \in \partial X, x,y \in X$.
The following Lemma is elementary:

\medskip

\begin{lemma} \label{busemann} We have $|g_{x,y}(\xi)| \leq
d_X(x,y)$ for all $\xi \in \partial x, x,y \in X$. Moreover
$g_{x,y}(\xi) = d_X(x,y)$ iff $y$ lies on the geodesic ray $[x,
\xi)$ while $g_{x,y}(\xi) = -d_X(x,y)$ iff $x$ lies on the
geodesic ray $[y, \xi)$.
\end{lemma}

\medskip

We recall the following Lemma from \cite{bourdon1}:

\medskip

\begin{lemma} \label{visualmvt} For $x, y \in X, \xi, \xi' \in
\partial X$ we have
$$
\rho_y(\xi, \xi') = \rho_x(\xi, \xi') f_{x,y}(\xi)^{1/2}
f_{x,y}(\xi')^{1/2}
$$
\end{lemma}

\medskip

An immediate corollary of the above Lemma is the following:

\medskip

\begin{lemma} \label{visualmoebius} The visual metrics $\rho_x, x \in X$
are Moebius equivalent to each other and
$$
\frac{d\rho_y}{d\rho_x} = f_{x,y}
$$
Hence the functions $f_{x,y}, g_{x,y}$ are continuous.
\end{lemma}

\medskip

It follows that the metric cross-ratio $[\xi\xi'\eta\eta']_{\rho_x}$ of a quadruple
$(\xi, \xi',\eta,\eta')$ is independent of the choice of $x \in
X$. Denoting this common value by $[\xi\xi'\eta\eta']$, it is
shown in \cite{bourdon2} that the cross-ratio is given by
$$
[\xi\xi'\eta\eta'] = \lim_{(a,a',b,b') \to (\xi, \xi',\eta,\eta')} \exp(\frac{1}{2}(d(a,b)+d(a',b') -
d(a,b') - d(a',b)))
$$
where the points $a,a',b,b' \in X$ converge radially towards
$\xi,\xi',\eta,\eta' \in \partial X$.

\medskip

We assume henceforth that $X$ is a proper, geodesically complete CAT(-1) space.
We let $\mathcal{M} = \mathcal{M}(\partial X, \rho_x)$ (this space is independent of the
choice of $x \in X$).

\medskip

\begin{lemma} The map
\begin{align*}
i_X : X & \to \mathcal{M} \\
         x & \mapsto \rho_x  \\
\end{align*}
is an isometric embedding and the image is closed in $\mathcal{M}$.
\end{lemma}

\medskip

\noindent{\bf Proof:} Given $x,y \in X$, extend $[x,y]$ to a
geodesic ray $[x, \xi)$ where $\xi \in \partial X$, then
$g_{x,y}(\xi) = d_X(x,y)$ hence $d_{\mathcal{M}}(\rho_x, \rho_y) =
\max_{\eta \in \partial X} g_{x,y}(\eta) = d_X(x,y)$, so $i_X$ is an isometric
embedding. Given $x_n \in X$ such that
$\rho_{x_n} \to \rho \in \mathcal{M}$, since $i_X$ is an isometry and the
sequence $\rho_{x_n}$ is bounded in $\mathcal{M}$,
so is the sequence $x_n$ in $X$. Passing to a subsequence we may assume $x_n \to a$ in
$X$, then $d_{\mathcal{M}}(\rho_{x_n}, \rho_a) = d_X(x_n,a) \to 0$ hence $\rho_a = \rho$. $\diamond$

\medskip

\subsection{Limiting comparison angles and derivatives of visual metrics}

\medskip

For points $a,x,a' \in X$ we denote by $\angle^{(-1)}axa' \in [0,\pi]$ the angle at the vertex corresponding to $x$ in a comparison triangle
  in $\mathbb{H}^2$ corresponding to the triangle $axa'$ in $X$. It is easy to show (see \cite{bourdon1}) that the map
  $X \times X \times X \to [0,\pi], (a,x,a') \mapsto \angle^{(-1)}axa'$ extends to a continuous map $\overline{X} \times X \times \overline{X} \to [0,\pi]$,
  so for $\xi,\xi' \in \partial X$ and $x \in X$ the limiting comparison angle $\angle^{(-1)}\xi x \xi'$ is defined, and
moreover
$$
\rho_x(\xi, \xi') = \sin(\angle^{(-1)}\xi x \xi'/2)
$$
For any point $y$ on the geodesic ray $[x,\xi)$ it follows easily from the
CAT(-1) inequality that
$$
\angle^{(-1)}y x \xi' \leq \angle^{(-1)}\xi x \xi'
$$
We note also that if a geodesic segment $[x,y]$ of length $\delta$
is common to both rays $[x,\xi)$ and $[x, \xi')$ then $\angle^{(-1)}y x \xi' = 0$ for $d(x,y) \leq \delta$.

\medskip

\begin{lemma} \label{visualderivative} For $x,y \in X$ and $\xi \in \partial X$, we have
$$
f_{x,y}(\xi) = \frac{1}{(e^t - e^{-t})\sin^2(\angle^{(-1)}y x \xi/2)+e^{-t}}
$$
\end{lemma}

\medskip

\noindent{\bf Proof:} Let $a$ tend to $\xi$ radially, let $r =
d_X(x,a), s = d_X(a, y)$ and let $\theta$ be the comparison angle
$\angle^{(-1)}y x a$. By the
hyperbolic law of cosine we have
$$
\cosh s = \cosh r \cosh t - \sinh r \sinh t \cos \theta
$$
which gives
$$
e^{s - r} + e^{-s - r} = (1 + e^{-2r})\frac{1}{2}(e^{t}+e^{-t}) -
\frac{1}{2}(1 - e^{-2r})(e^t - e^{-t})\cos \theta
$$
Now as $r \to \infty$ we have $s \to \infty$, and by definition $r
- s \to B(\xi, x, y)$, also $\theta \to \angle^{(-1)}y x \xi$,
hence letting $r \to \infty$ above gives
$$
\frac{1}{f_{x,y}(\xi')} = \frac{1}{2}(e^t + e^{-t}) -
\frac{1}{2}(e^t - e^{-t})\cos(\angle^{(-1)}y x \xi) = (e^t - e^{-t})\sin^2(\angle^{(-1)}y x \xi)/2)+e^{-t}
$$
$\diamond$

\medskip

We now consider the behaviour of the derivatives $f_{x,y}$ as $t = d(x,y) \to 0$ and the point $y$
converges radially towards $x$. For functions $F_t$ on $\partial X$ we write $F_t = o(t)$ if
$||F_t||_{\infty} = o(t)$. We have the following formula, which may be thought of as a formula
for the derivative of the map $i_X$ along a geodesic:

\medskip

\begin{lemma} \label{embedderiv} As $t \to 0$ we have
$$
g_{x,y}(\xi) = t \cos(\angle^{(-1)}y x \xi) + o(t)
$$
\end{lemma}

\medskip

\noindent{\bf Proof:} As $t \to 0$ we have
\begin{align*}
g_{x,y}(\xi) & = -\log((e^t - e^{-t})\sin^2(\angle^{(-1)}y x \xi/2)+e^{-t}) \\
& = -\log(2t\sin^2(\angle^{(-1)}y x \xi/2) + 1 - t +
o(t)) \\
& = -(2t\sin^2(\angle^{(-1)}y x \xi/2) - t) + o(t) \\
& = t \cos(\angle^{(-1)}y x \xi) + o(t) \\
\end{align*}
$\diamond$

\medskip

\section{Geodesic conjugacies, Moebius maps, conformal maps, and the integrated
Schwarzian}

\medskip

We start by recalling the definitions of conformal maps, Moebius
maps, and the abstract geodesic flow of a CAT(-1) space.

\medskip

\begin{definition} A homeomorphism between metric spaces $f :
(Z_1, \rho_1) \to (Z_2, \rho_2)$ with no isolated points is said to be {\it conformal} if
for all $\xi \in Z_1$, the limit
$$
df_{\rho_1, \rho_2}(\xi) := \lim_{\eta \to \xi} \frac{\rho_2(f(\xi),
f(\eta))}{\rho_1(\xi, \eta)}
$$
exists and is positive. The positive function $df_{\rho_1,
\rho_2}$ is called the derivative of $f$ with respect to $\rho_1, \rho_2$.
We say $f$ is {\it $C^1$ conformal} if its derivative is continuous.

\medskip

Two metrics $\rho_1, \rho_2$ inducing the same topology on a set
$Z$, such that $Z$ has no isolated points,
are said to be conformal (respectively $C^1$ conformal) if the
map $id_Z : (Z, \rho_1) \to (Z, \rho_2)$ is conformal
(respectively $C^1$ conformal). In this case we denote the
derivative of the identity map by $\frac{d\rho_2}{d\rho_1}$.
\end{definition}

\medskip

\begin{definition} A homeomorphism between metric spaces $f :
(Z_1, \rho_1) \to (Z_2, \rho_2)$ (where $Z_1$ has at least four
points) is said to be Moebius if it preserves metric cross-ratios
with respect to $\rho_1, \rho_2$. The derivative of $f$ is defined
to be the derivative $\frac{df_*\rho_2}{\rho_1}$ of the Moebius
equivalent metrics $f_* \rho_2, \rho_1$ as defined in section 2
(where $f_* \rho_2$ is the pull-back of $\rho_2$ under $f$).
\end{definition}

\medskip

From the results of section 2 it follows that any Moebius map
between compact metric spaces with no isolated points is $C^1$ conformal,
and the two definitions of the derivative of $f$ given above
coincide. Moreover any Moebius map $f$ satisfies the geometric
mean-value theorem,
$$
\rho_2(f(\xi), f(\eta))^2 = \rho_1(\xi,\eta)^2
df_{\rho_1,\rho_2}(\xi) df_{\rho_1,\rho_2}(\xi)
$$

\medskip

\begin{definition} Let $(X, d)$ be a CAT(-1) space. The abstract geodesic flow
space of $X$ is defined to be the space of bi-infinite geodesics
in $X$,
$$
\mathcal{G}X := \{ \gamma : (-\infty,+\infty) \to X | \gamma
\hbox{ is an isometric embedding} \}
$$
endowed with the topology of uniform convergence on compact
subsets. This topology is metrizable with a distance defined by
$$
d_{\mathcal{G}X}(\gamma_1, \gamma_2):= \int_{-\infty}^{\infty}
d(\gamma_1(t), \gamma_2(t)) \frac{e^{-|t|}}{2} \ dt
$$
We define also a projection
\begin{align*}
\pi_X : \mathcal{G}X & \to X \\
             \gamma  & \mapsto \gamma(0) \\
\end{align*}
It is shown in Bourdon \cite{bourdon1} that $\pi_X$ is
$1$-Lipschitz.

\medskip

The abstract geodesic flow of $X$ is defined to be the one-parameter group
of homeomorphisms
\begin{align*}
\phi^X_t : \mathcal{G}X & \to \mathcal{G}X \\
         \gamma       & \mapsto \gamma_t \\
\end{align*}
for $t \in \mathbb{R}$, where $\gamma_t$ is the geodesic $s
\mapsto \gamma(s+t)$.

\medskip

The flip is defined to be the map
\begin{align*}
\mathcal{F}_X : \mathcal{G}X & \to \mathcal{G}X \\
              \gamma & \mapsto \overline{\gamma} \\
\end{align*}
where $\overline{\gamma}$ is the geodesic $s
\mapsto \gamma(-s)$.
\end{definition}

\medskip

We observe that for a simply connected complete Riemannian manifold $X$ with
sectional curvatures bounded above by $-1$, the map
\begin{align*}
\mathcal{G}X & \to T^1 X \\
 \gamma      & \mapsto \gamma'(0) \\
\end{align*}
is a homeomorphism conjugating the abstract geodesic flow of $X$
to the usual geodesic flow of $X$ and the flip $\mathcal{F}$ to the
usual flip on $T^1 X$.

\medskip

We note that that for any CAT(-1) space $X$ there is a continuous
surjection
\begin{align*}
\mathcal{E}_X : \mathcal{G}X & \to \partial^2 X \\
                   \gamma  & \mapsto (\gamma(-\infty),
                   \gamma(+\infty)) \\
\end{align*}
which induces a homeomorphism $\mathcal{G}X / (\phi_t)_{t \in
\mathbb{R}} \to \partial^2 X$. Following Bourdon \cite{bourdon1}, we have the following:

\medskip

\begin{prop} \label{confconj} Let $f : \partial X \to \partial Y$ be a
conformal map between the boundaries of CAT(-1) spaces $X, Y$
equipped with visual metrics. Then $f$ induces a bijection $\phi_f :
\mathcal{G}X \to \mathcal{G}Y$ conjugating the geodesic flows,
which is a homeomorphism if $f$ is $C^1$ conformal. If $f$ is
Moebius then $\phi_f$ is flip-equivariant.
\end{prop}

\medskip

\noindent {\bf Proof:} Given $\gamma \in \mathcal{G}X$, let
$\mathcal{E}_X(\gamma) = (\xi, \eta), x = \gamma(0)$, then there
is a unique point $y \in (f(\xi),f(\eta))$ such that $df_{\rho_x,
\rho_y}(\eta) = 1$. Define $\phi_f(\gamma) = \gamma^*$ where
$\gamma^*$ is the unique geodesic in $Y$ satisfying $\mathcal{E}_Y(\gamma^*) = (f(\xi),
f(\eta)), \gamma^*(0) = y$. Then $\phi_f : \mathcal{G}X \to
\mathcal{G}Y$ is a bijection conjugating the geodesic flows.

\medskip

\noindent{\bf Claim.} The map $\phi_f$ is continuous if $f$ is $C^1$
conformal.

\medskip

\noindent{\bf Proof of Claim:} Let $\gamma_n \to \gamma$ in
$\mathcal{G}X$. Let $x = \gamma(0), x_n = \gamma_n(0),
\mathcal{E}_X(\gamma) = (\xi,\eta), \mathcal{E}_X(\gamma_n) =
(\xi_n,\eta_n)$. Then $x_n \to x, (\xi_n,\eta_n) \to (\xi, \eta)$,
hence
\begin{align*}
\rho_x(\xi_n,\eta_n) & = \rho_{x_n}(\xi_n,\eta_n)
\frac{d\rho_x}{d\rho_{x_n}}(\xi_n)^{1/2}
\frac{d\rho_x}{d\rho_{x_n}}(\eta_n)^{1/2} \\
                     & = \frac{d\rho_x}{d\rho_{x_n}}(\xi_n)^{1/2}
\frac{d\rho_x}{d\rho_{x_n}}(\eta_n)^{1/2} \\
                     & \to 1 \\
\end{align*}
since $|g_{x_n,x}| \leq d(x,x_n) \to 0$. Letting $y = \pi_Y \circ \phi_f(\gamma)$, this implies
$\rho_y(f(\xi_n),f(\eta_n)) \to \rho_y(f(\xi),f(\eta)) = 1$ since $f$
is continuous.

\medskip

Fix $\epsilon > 0$ small and $n$ large such that $\rho_y(f(\xi_n), f(\eta_n)) \geq 1 - \epsilon$.
If $a_t, b_t$ are points converging radially towards $f(\xi_n),
f(\eta_n)$, then as $t \to +\infty$ there are points
$\overline{z_t}$ in the comparison triangle
$\overline{a_t}\overline{y}\overline{b_t}$ on the side
$\overline{a_t}\overline{b_t}$ such that $d(\overline{z_t},
\overline{y}) \leq C(\epsilon)$ for some constant $C(\epsilon)$
which tends to $0$ as $\epsilon$ tends to $0$. Hence we obtain a
point $z_n \in (f(\xi_n), f(\eta_n))$ such that $d(z_n, y) \leq
C(\epsilon)$. Therefore $d(z_n, y) \to 0$ as $n \to \infty$.

\medskip

Let $z^*_n = \pi_Y \circ \phi_f(\gamma_n)$. Then since $z_n, z^*_n$
both lie on the geodesic $\phi_f(\gamma_n)$ and
$df_{\rho_{x_n},\rho_{z^*_n}}(\eta_n) = 1$, we have
\begin{align*}
d(z^*_n, z_n) & = |\log df_{\rho_{x_n}, \rho_{z_n}}(\eta_n)| \\
              & = \left|\log \left( df_{\rho_x, \rho_y}(\eta_n) \frac{d\rho_x}{d\rho_{x_n}}(\eta_n)
              \frac{d\rho_{z_n}}{d\rho_{y}}(f(\eta_n)) \right) \right| \\
              & \to |\log(1 \cdot 1 \cdot 1)| = 0 \\
\end{align*}
since $f$ is $C^1$ conformal with $df_{\rho_x, \rho_y}(\eta) = 1$ and $\eta_n \to
\eta$, $d(x_n, x)+d(z_n,y) \to 0$. Hence the basepoints $z^*_n$ of the geodesics
$\phi_f(\gamma_n)$ converge to the basepoint $y$ of the geodesic $\phi_f(\gamma)$, and
the endpoints $(f(\xi_n), f(\eta_n))$ of $\phi_f(\gamma_n)$ converge
to the endpoints $(f(\xi), f(\eta))$ of $\phi_f(\gamma)$, from which it follows easily
that $\phi_f(\gamma_n) \to \phi_f(\gamma)$ in $\mathcal{G}Y$. This
finishes the proof of the Claim.

\medskip

Since the inverse of a $C^1$ conformal map is clearly $C^1$
conformal, $f^{-1}$ also induces a continuous conjugacy $\psi_f :
\mathcal{G}Y \to \mathcal{G}X$ which is clearly inverse to $\phi_f$,
hence $\phi_f$ is a homeomorphism if $f$ is $C^1$ conformal.

\medskip

If $f$ is Moebius, then with the same notation as
above, by the geometric mean-value theorem we have $df_{\rho_x, \rho_y}(\xi) df_{\rho_x, \rho_y}(\eta) = 1$,
hence $df_{\rho_x, \rho_y}(\xi)= 1$, and it follows that $\phi_f$ is
flip-equivariant. $\diamond$

\medskip

The proof of flip-equivariance of the conjugacy for a Moebius map above
motivates the following definition:

\medskip

\begin{definition} Let $f : \partial X \to \partial Y$ be a
conformal map between boundaries of CAT(-1) spaces equipped with
visual metrics. The integrated Schwarzian of $f$ is the function
$S(f) : \partial^2 X \to \mathbb{R}$ defined by
$$
S(f)(\xi, \eta) := - \log (df_{\rho_x, \rho_y}(\xi) df_{\rho_x,
\rho_y}(\eta)) \ \ (\xi,\eta) \in \partial^2 X
$$
where $x,y$ are any two points $x \in (\xi, \eta), y \in (f(\xi),
f(\eta))$ (it is easy to see that the quantity defined above is
independent of the choices of $x$ and $y$).
\end{definition}

\medskip

We note that $S(f)$ is continuous if $f$ is $C^1$ conformal, and
for any $\gamma \in \mathcal{G}X$ with $\mathcal{E}_X(\gamma) =
(\xi, \eta)$, we have
$$
\phi_f(\mathcal{F}_X(\gamma)) =
\mathcal{F}_Y(\phi^Y_{-t}(\phi_f(\gamma)))
$$
where $t = S(f)(\xi, \eta)$, hence the integrated Schwarzian of $f$
measures the deviation of the induced conjugacy $\phi_f$ from being flip-equivariant.

\medskip

We consider now the relation between the integrated Schwarzian and the continuity of the conjugacy
$\phi_f$ near infinity. In particular we consider the continuity properties of $\phi_f$ along geodesics.

\medskip

\begin{definition} Let $X$ be a simply connected complete Riemannian manifold with sectional curvatures bounded above and
below, $-b^2 \leq K \leq -1$. A sequence of pairs of unit tangent vectors $(v_n, w_n) \in T^1 X \times T^1 X$ is said to be
forward asymptotic along a geodesic $\gamma \in \mathcal{G}X$ if:

\medskip

\noindent 1) There are times $t_n \to +\infty$ such that
$v_n = \gamma'(t_n)$ and $d_{T^1 X}(v_n, w_n)  \to 0$
(the distance on $T^1 X$ being the Sasaki metric).

\medskip

\noindent 2) Let $\gamma_n \in \mathcal{G}X$ such that $\gamma'_n(0) = w_n$, let $\mathcal{E}_X(\gamma) = (\xi, \eta),
\mathcal{E}_X(\gamma_n) = (\xi_n, \eta_n)$. Then we require $\xi_n \to \xi_0 \neq \eta$ as $n \to \infty$.
\end{definition}

\medskip

We have:

\medskip

\begin{prop} \label{fwdasymptotic} Let $X$ be a simply connected complete Riemannian manifold
with sectional curvatures bounded above and below, $-b^2 \leq K \leq -1$, and let $Y$ be a CAT(-1) space. Let
$f : \partial X \to \partial Y$ be a $C^1$ conformal map and $\phi = \phi_f : T^1 X \to \mathcal{G}Y$ the associated
geodesic conjugacy. Then for any sequence $(v_n, w_n)$ forward asymptotic
along a geodesic $\gamma$, we have
$$
d_Y(\pi_Y \circ \phi(v_n), \pi_Y \circ \phi(w_n)) \to 0
$$
\end{prop}

\medskip

\noindent{\bf Proof:}
Let  $(v_n, w_n)$ be a forward asymptotic sequence along a geodesic $\gamma$, so there are times $t_n \to +\infty$ such that
$v_n = \gamma'(t_n)$ and $d_{T^1 X}(v_n, w_n) \to 0$. Let
$x = \gamma(0)$, $x_n = \gamma(t_n) \in X$, $y = \pi_Y \circ \phi(\gamma'(0))$,
$y_n = \pi_Y \circ \phi(v_n) \in Y$.
Let $\gamma_n \in \mathcal{G}X$ with $\gamma'_n(0) = w_n$, let $\mathcal{E}_X(\gamma_n) = (\xi_n,
\eta_n), \mathcal{E}_X(\gamma) = (\xi, \eta)$, then by hypothesis
$\xi_n \to \xi_0 \neq \eta$. Since the
curvature of $X$ is bounded below by $-b^2$, for any $T \in
\mathbb{R}$ the time-$T$-map of the geodesic flow $\phi^X_T : T^1
X \to T^1 X$ is Lipschitz. This follows from the fact
that the differential of the map $\phi^X_T$ is given in terms of Jacobi
fields and their derivatives, and by well known comparison
arguments, Jacobi fields in $X$ grow at most as fast as Jacobi
fields in the hyperbolic space of constant curvature $-b^2$, hence
 $||d\phi^X_T||$ is bounded on $T^1 X$. It follows that for any fixed large $T$,
 $d_{T^1 X}(\phi^X_T(v_n), \phi^X_T(w_n)) \to 0$, hence the visual distance
  $\rho_{x_n}(\eta, \eta_n) \to 0$. It is easy to see that this also implies $\rho_x(\eta, \eta_n) \to 0$.

\medskip

\noindent{\bf Claim.} We have
$$
\lim_{n \to \infty} \frac{d\rho_{x_n}}{d\rho_x}(\eta_n) e^{-t_n} =
\lim_{n \to \infty} \frac{d\rho_{y_n}}{d\rho_y}(f(\eta_n))
e^{-t_n} = 1
$$

\noindent{\bf Proof of Claim:} Fix $\epsilon > 0$ small. Let $\alpha_n \in \mathcal{G}X$ be
a geodesic with $\alpha_n(0) = x_n, \alpha_n(+\infty) = \eta_n$.
Then the Riemannian angle between $\alpha'_n(0), v_n$
tends to $0$ (since the comparison angle $\angle^{(-1)} \eta_n x_n \eta$ tends to
$0$), so the Riemannian angle between $\alpha'_n(0),
-v_n$ tends to $\pi$. Hence the limit of comparison angles
$(\lim_{t \to +\infty} \angle \overline{\alpha_n(t)}\overline{x_n}\overline{x})$ tends to $\pi$ as
$n \to \infty$ (where $\overline{\alpha_n(t)}\overline{x_n}\overline{x}$ is a comparison triangle in $\mathbb{H}^2$).
Fix $n$ large such that this limiting angle is
larger than $\pi - \epsilon$. For $t > 0$ large the comparison triangles
$\overline{\alpha_n(t)}\overline{x_n}\overline{x}$ in $\mathbb{H}^2$ have an
angle at the vertex $\overline{x_n}$ greater than $\pi -
\epsilon$, hence the sides satisfy
$$
d(\alpha_n(t), x) - d(\alpha_n(t), x_n) \geq d(x_n,x) - C(\epsilon)
$$
for some constant $C(\epsilon)$ which tends to $0$ as $\epsilon$ tends to
$0$. Letting $t \to +\infty$, we have $B(\eta_n, x, x_n) \geq t_n
- C(\epsilon)$, hence
$$
e^{t_n} = e^{d(x,x_n)} \geq \frac{d\rho_{x_n}}{d\rho_x}(\eta_n)
\geq e^{-C(\epsilon)}e^{t_n}
$$
therefore $\frac{d\rho_{x_n}}{d\rho_x}(\eta_n) e^{-t_n} \to 1$.

\medskip

Now using the geometric mean value theorem for visual metrics we
have
\begin{align*}
\rho_{y_n}(f(\eta_n),f(\eta)) & = \frac{\rho_{y_n}(f(\eta_n),
f(\eta))}{\rho_y(f(\eta_n), f(\eta))} \frac{\rho_y(f(\eta_n),
f(\eta))}{\rho_x(\eta_n,\eta)}
\frac{\rho_x(\eta_n,\eta)}{\rho_{x_n}(\eta_n,\eta)}
\rho_{x_n}(\eta_n,\eta) \\
& = \left(e^{t_n}
\frac{d\rho_{y_n}}{d\rho_y}(\eta_n)\right)^{1/2} \frac{\rho_y(f(\eta_n),
f(\eta))}{\rho_x(\eta_n,\eta)} \left(e^{-t_n}
\left(\frac{d\rho_{x_n}}{d\rho_x}(\eta_n)\right)^{-1}\right)^{1/2}
\rho_{x_n}(\eta_n, \eta) \\
& \leq \frac{\rho_y(f(\eta_n),
f(\eta))}{\rho_x(\eta_n,\eta)} \left(e^{t_n} \left(\frac{d\rho_{x_n}}{d\rho_x}(\eta_n)\right)^{-1}\right)^{1/2}
\rho_{x_n}(\eta_n, \eta) \\
& \to 1 \cdot 1 \cdot 0 = 0 \\
\end{align*}

\medskip

Now $\rho_{y_n}(f(\eta_n), f(\eta)) \to 0$ and $d(y_n,y) = t_n$
implies that
$$
\lim_{n \to \infty} \frac{d\rho_{y_n}}{d\rho_y}(f(\eta_n))
e^{-t_n} = 1
$$
by the same argument used above to show that $\frac{d\rho_{x_n}}{d\rho_x}(\eta_n) e^{-t_n} \to
1$. This finishes the proof of the Claim.

\medskip

Now note that since $f(\xi_n) \to f(\xi_0) \neq f(\eta)$ and $y_n
\to \eta$ radially, we have $\rho_{y_n}(f(\xi), f(\xi_n)) \to 0$.
Hence
\begin{align*}
\rho_{y_n}(f(\xi_n), f(\eta_n)) & \geq \rho_{y_n}(f(\xi), f(\eta)) -
\rho_{y_n}(f(\xi), f(\xi_n)) - \rho_{y_n}(f(\eta), f(\eta_n)) \\
      & = 1 - \rho_{y_n}(f(\xi), f(\xi_n)) - \rho_{y_n}(f(\eta), f(\eta_n)) \\
      & \to 1 \\
\end{align*}.
Fix $\epsilon > 0$ small. Fix $n$ large such that $\rho_{y_n}(f(\xi_n), f(\eta_n)) \geq 1 - \epsilon$.
If $a_t, b_t$ are points converging radially towards $f(\xi_n),
f(\eta_n)$, then as $t \to +\infty$ there are points
$\overline{z_t}$ in the comparison triangle
$\overline{a_t}\overline{y_n}\overline{b_t}$ on the side
$\overline{a_t}\overline{b_t}$ such that $d(\overline{z_t},
\overline{y_n}) \leq C(\epsilon)$ for some constant $C(\epsilon)$
which tends to $0$ as $\epsilon$ tends to $0$. Hence we obtain a
point $z_n \in (f(\xi_n), f(\eta_n))$ such that $d(z_n, y_n) \leq
C(\epsilon)$. Therefore $d(z_n, y_n) \to 0$ as $n \to \infty$.

\medskip

Let $x^*_n = \pi_X(w_n), z^*_n = \pi_Y \circ \phi(w_n)$. Note $d(x^*_n, x_n) \to 0$.
Since $z_n, z^*_n$ lie on the geodesic $(f(\xi_n), f(\eta_n))$ and $df_{\rho_{x^*_n}, \rho_{z^*_n}}(\eta_n) = 1$,
we have
\begin{align*}
d(z^*_n, z_n) & = \left|\log df_{\rho_{x^*_n}, \rho_{z_n}}(\eta_n) \right| \\
              & = \left|\log \left( df_{\rho_x, \rho_y}(\eta_n) \left( \frac{d\rho_x}{d\rho_{x_n}}(\eta_n)
              \frac{d\rho_{y_n}}{d\rho_y}(f(\eta_n)) \right) \left( \frac{d\rho_{x_n}}{d\rho_{x^*_n}}(\eta_n)
              \frac{d\rho_{z_n}}{d\rho_{y_n}}(f(\eta_n)) \right)            \right) \right| \\
              & \to |\log(1 \cdot 1 \cdot 1)| = 0 \\
\end{align*}
since $f$ is $C^1$ conformal with $df_{\rho_x, \rho_y}(\eta) = 1$ and $\eta_n \to
\eta$, $d(x^*_n, x_n)+d(z_n,y_n) \to 0$, and the term in the middle of the product tends to $1$ by the Claim
proved earlier.

\medskip

Hence $d(\pi_Y \circ \phi(v_n), \pi_Y \circ \phi(w_n)) = d(y_n,z^*_n)
\to 0$. $\diamond$

\medskip

\begin{prop} \label{distancediff} Let $X,Y, f, \phi$ be as in the previous Proposition.
Let $x \in X$ and $(\xi,\eta) \in \partial^2 X$. Let $\alpha, \beta : [0,\infty) \to X$
be geodesic rays joining $x$ to $\xi, \eta$ respectively. Let $x_t = \alpha(t), y_t = \beta(t), v_t = \alpha'(t), w_t = \beta'(t)$, then
$$
d_Y(\pi_Y \circ \phi(v_t), \pi_Y \circ \phi(w_t)) - d_X(x_t, y_t) \to S(f)(\xi, \eta)
$$
as $t \to +\infty$.
\end{prop}

\medskip

\noindent{\bf Proof:} Let $\gamma_t$ be the bi-infinite geodesic passing through $x_t, y_t$, with endpoints $(\xi_t, \eta_t) \in \partial^2 X$,
so that $(\xi_t, \eta_t) \to (\xi, \eta)$ as $t \to +\infty$. Let $v'_t, w'_t$ be the tangent vectors to $\gamma_t$ at the points $x_t, y_t$
pointing respectively towards $\xi_t, \eta_t$. Then it is a standard fact that for any sequence $t_n \to +\infty$,
the sequences of pairs $\{( v_{t_n}, v'_{t_n})\}, \{(w_{t_n}, w'_{t_n})\}$ are forward asymptotic along $\alpha, \beta$ respectively. Letting
$p_{t} = \pi_Y \circ \phi(v_{t}), q_{t} = \pi_Y \circ \phi(w_t), p'_t = \pi_Y \circ \phi(v'_t), q'_t = \pi_Y \circ \phi(w'_t)$, then by
Proposition \ref{fwdasymptotic} we have $d_Y(p_{t_n}, p'_{t_n}) \to 0, d_Y(q_{t_n}, q'_{t_n}) \to 0$ as $n \to \infty$. By definition of the integrated
Schwarzian, we have $d_Y(p'_{t_n}, q'_{t_n}) = d_X(x_{t_n}, y_{t_n}) + S(f)(\xi_{t_n}, \eta_{t_n})$, since $S(f)$ is continuous it follows that
$d_Y(p'_{t_n}, q'_{t_n}) - d_X(x_{t_n}, y_{t_n}) = S(f)(\xi_{t_n}, \eta_{t_n}) \to S(f)(\xi, \eta)$ as $n \to \infty$. The result follows. $\diamond$

\medskip

We can now prove Theorem \ref{crossratiodistn}:

\medskip

\noindent{\bf Proof:} We first note that $f : U \to V$ induces a geodesic conjugacy between the flow invariant subsets of $\mathcal{G}X, \mathcal{G}Y$
with endpoints in $U,V$ respectively, for which the same arguments as above show that the conclusion of Proposition \ref{distancediff}
above holds. Fix a basepoint $x \in X$.
Now given $(\xi, \xi', \eta, \eta') \in \partial^4 U$, let $\alpha, \beta, \gamma, \delta$ be geodesic rays joining $x$ to $\xi, \eta, \xi', \eta'$
respectively. Let $x_t = \alpha(t), y_t = \beta(t), a_t = \gamma(t), b_t = \delta(t)$, let $v_t = \alpha'(t), w_t = \beta'(t), v'_t = \gamma'(t),
w'_t = \delta'(t)$ and let $p_t = \pi_Y \circ \phi(v_t), q_t = \pi_Y \circ \phi(w_t), r_t = \pi_Y \circ \phi(v'_t), s_t = \pi_Y \circ \phi(w'_t)$.
Then the points $p_t, q_t, r_t, s_t$ converge radially towards $f(\xi), f(\eta), f(\xi'), f(\eta')$, hence

\begin{align*}
\log \frac{[f(\xi), f(\xi'), f(\eta), f(\eta')]}{[\xi, \xi',\eta,\eta']}  = & \frac{1}{2} (\lim_{t \to \infty}
(d_Y(p_t, q_t) - d_X(x_t,y_t)) + (d_Y(r_t, s_t) - d_X(a_t, b_t)) \\
& - (d_Y(p_t, s_t) - d_X(x_t, b_t)) - (d_Y(q_t, r_t) - d_X(y_t, a_t)) ) \\
= & \frac{1}{2}(S(f)(\xi,\eta) + S(f)(\xi', \eta') - S(f)(\xi, \eta') - S(f)(\xi', \eta)) \\
\end{align*}

(using Proposition \ref{distancediff} in the last line above) $\diamond$.

\medskip

\begin{definition} Let $X$ be a simply connected complete Riemannian manifold
with sectional curvatures bounded above and below, $-b^2 \leq K \leq -1$, and let $Y$ be a CAT(-1) space. A homeomorphism
$\phi : T^1 X \to \mathcal{G}Y$ is said to be uniformly
continuous along geodesics, if, given $\gamma \in \mathcal{G}X$,
and a sequence $(v_n, w_n) \in T^1 X \times T^1 X$ which is forward asymptotic along $\gamma$, we have
$$
d(\pi_Y \circ \phi(v_n), \pi_Y \circ \phi(w_n)) +
d(\pi_Y \circ \phi(-v_n),
\pi_Y \circ \phi(-w_n)) \to 0
$$
\end{definition}

\medskip

We note that any uniformly continuous homeomorphism $\phi : T^1 X
\to \mathcal{G}Y$ is uniformly continuous along geodesics. We can
now prove Theorem \ref{geodesicmoebius}:

\medskip

\noindent {\bf Proof of Theorem \ref{geodesicmoebius}}: We first
note that if $\gamma_1, \gamma_2 \in \mathcal{G}X$ are geodesics
with $\gamma_1(+\infty) = \gamma_2(+\infty)$, then it follows easily from
the definition of uniform continuity along geodesics that
$\phi(\gamma'_1(0))(+\infty) = \phi(\gamma'_2(0))(+\infty)$. Hence there is a map $f : \partial X \to \partial Y$
such that $\mathcal{E}_Y(\phi(v)) = (f(\xi), f(\eta))$ where
$(\xi, \eta) = \mathcal{E}_X(\gamma)$, $\gamma \in
\mathcal{G}X$ being such that $\gamma'(0) = v$, and it is not hard to show that
$f$ is continuous. Moreover $f$ is surjective since
$Y$ is geodesically complete and $\phi$ is surjective. Also given $(\xi,\eta) \in \partial^2
X$, choosing $\gamma$ with $\mathcal{E}_X(\gamma) = (\xi,\eta)$,
we have $(f(\xi), f(\eta)) =  \mathcal{E}_Y(\phi(\gamma'(0))) \in \partial^2 Y$,
in particular $f(\xi) \neq f(\eta)$.
Thus $f$ is injective, and since $\partial X, \partial Y$ are
compact Hausdorff spaces, $f$ is a homeomorphism.

\medskip

Given a
quadruple of distinct points $(\xi, \xi', \eta, \eta') \in
\partial^4 X$, let $\gamma_1, \gamma_2$ be geodesics with
$\mathcal{E}_X(\gamma_1) = (\xi, \eta), \mathcal{E}_X(\gamma_2) = (\xi',
\eta')$, and let $t_n \to +\infty$. Let $a_n = \gamma_1(-t_n),
a'_n = \gamma_2(-t_n), b_n = \gamma_1(t_n), b'_n = \gamma_2(t_n)$
so
$$
[\xi\xi'\eta\eta'] = \lim_{n \to \infty}
\exp(\frac{1}{2}(d(a_n,b_n)+d(a'_n,b'_n)-d(a_n,b'_n)-d(a'_n,b_n)))
$$
Let $\alpha_n = \pi_Y \circ \phi (\gamma'_1(-t_n)), \alpha'_n =
\pi_Y \circ \phi (\gamma'_2(-t_n)), \beta_n = \pi_Y \circ \phi
(\gamma'_1(t_n)), \beta'_n = \pi_Y \circ \phi (\gamma'_2(t_n))$,
so that
$$
[f(\xi)f(\xi')f(\eta)f(\eta')] = \lim_{n \to \infty}
\exp(\frac{1}{2}(d(\alpha_n,\beta_n)+d(\alpha'_n,\beta'_n)-d(\alpha_n,\beta'_n)-d(\alpha'_n,\beta_n)))
$$
Note that $d(a_n, b_n) = d(\alpha_n, \beta_n), d(a'_n, b'_n) = d(\alpha'_n,
\beta'_n)$ since $\phi$ conjugates the geodesic flows. Clearly
the Theorem follows from the following claim:

\medskip

\noindent{\bf Claim.} We have $d(a_n, b'_n) - d(\alpha_n,
\beta'_n) \to 0, d(a'_n, b_n) - d(\alpha'_n,
\beta_n) \to 0$ as $n \to \infty$.

\medskip

\noindent{\bf Proof of Claim.} Let $\gamma_n : [0,l_n] \to X$ be
the geodesic segment with $\gamma_n(0) = a_n, \gamma_n(l_n) =
b'_n$, where $l_n = d(a_n, b'_n)$. Then it is a standard fact that
the Riemannian angle between the vectors $\gamma'_1(-t_n), v_n = \gamma'_n(0)$
tends to $0$, as does the angle between the vectors
$\gamma'_2(t_n), w_n = \gamma'_n(l_n)$. Letting $p_n = \pi_Y \circ
\phi(v_n), q_n = \pi_Y \circ \phi(w_n)$, we have $d(p_n, q_n) =
d(a_n, b'_n)$ since $\phi$ is a geodesic conjugacy. Moreover since
$\phi$ is uniformly continuous along geodesics, it follows that
$d(p_n, \alpha_n) \to 0, d(q_n, \beta'_n) \to 0$. Hence $d(a_n, b'_n) - d(\alpha_n,
\beta'_n) \to 0$ and a similar argument shows $d(a'_n, b_n) - d(\alpha'_n,
\beta_n) \to 0$. $\diamond$

\medskip

\noindent{\bf Proof of Theorem \ref{mobconfgeo}:} The forward implication follows from Theorem \ref{geodesicmoebius}.

\medskip

For the backward implication, given $f : \partial X \to \partial Y$ a Moebius map, let $\phi : T^1 X \to \mathcal{G}Y$ denote the
induced conjugacy of geodesic flows given by Proposition
\ref{confconj}. We show that $\phi$ is uniformly continuous along
geodesics:

\medskip

Let $(v_n, w_n)$ be a forward asymptotic sequence. By Proposition \ref{fwdasymptotic}, we have $d(\pi_Y \circ \phi(v_n), \pi_Y \circ \phi(w_n)) \to 0$.
Since $f$ is Moebius, the conjugacy $\phi$ is flip-equivariant, hence $\pi_Y \circ \phi(-v_n) = \pi_Y \circ \phi(v_n), \pi_Y \circ \phi(-w_n) = \pi_Y \circ \phi(w_n)$, thus $d(\pi_Y \circ \phi(-v_n), \pi_Y \circ \phi(-w_n)) = d(\pi_Y \circ \phi(v_n), \pi_Y \circ \phi(w_n)) \to 0$. $\diamond$

\medskip

It follows from the chain rule that the integrated Schwarzian
satisfies the following transformation rule: given conformal maps
$f : \partial X \to \partial Y, g: \partial Y \to \partial Z$, where $X, Y, Z$ are CAT(-1) spaces,
we have
$$
S(g \circ f) = S(g) \circ f + S(f)
$$
For the group $G$ of $C^1$ conformal self-maps of the boundary $\partial X$ of a CAT(-1) space,
the map
\begin{align*}
c : G & \to C(\partial^2 X) \\
    f & \mapsto S(f)          \\
\end{align*}
is therefore a $G$-cocycle with values in
the vector space $C(\partial^2 X)$ of continuous functions on
$\partial^2 X$ endowed with its natural $G$-action.

\medskip

For the group $G$ of $C^1$ conformal self-maps of $V \subset \partial Y$, it follows
that the subgroup $ker \ c := \{ g \in G | S(g) = 0 \} < G$ coincides with the group
of Moebius self-maps of $V$. Hence for conformal maps $f,g: U \to
V$, $g \circ f^{-1}$ is Moebius if and only if $S(g \circ f^{-1})
= 0$. Using the identities
\begin{align*}
S(g \circ f^{-1}) & = S(g) \circ f^{-1} + S(f^{-1}) \\
0 = S(f \circ f^{-1}) & = S(f) \circ f^{-1} + S(f^{-1}) \\
\end{align*}
it follows that $S(g \circ f^{-1}) = (S(g) - S(f)) \circ f^{-1}$,
hence $f,g$ differ by post-composition with a Moebius map if and
only if $S(g) = S(f)$.

\medskip

For $C^1$ conformal maps $f : \partial X \to \partial Y$ such that the integrated Schwarzian
$S(f)$ is bounded, and $X$ is a simply connected manifold with pinched negative sectional curvatures, we have
the following version of the geometric mean value theorem:

\medskip

\begin{theorem}\label{confgmvt} Let $X$ be a simply connected complete Riemannian manifold
with sectional curvatures bounded above and below, $-b^2 \leq K \leq -1$, and let $Y$ be a CAT(-1) space. Let
$f : \partial X \to \partial Y$ be a $C^1$ conformal map such that $S(f)$ is bounded. Then for all $(\xi, \eta) \in \partial^2 X$
and $x \in X, y \in Y$, we have
$$
e^{-4||S(f)||_{\infty}} df_{\rho_x,\rho_y}(\xi) df_{\rho_x,\rho_y}(\eta) \leq \left(\frac{\rho_y(f(\xi), f(\eta))}{\rho_x(\xi, \eta)}\right)^2
\leq e^{4||S(f)||_{\infty}} df_{\rho_x,\rho_y}(\xi) df_{\rho_x,\rho_y}(\eta)
$$
\end{theorem}

\medskip

\noindent{\bf Proof:} Fix $x \in X, y \in Y$. For a triple $(\xi, \xi', \eta') \in \partial^3 X$, we define
$$
\delta(\xi, \xi', \eta') := \frac{\rho_y(f(\xi), f(\xi'))\rho_y(f(\xi), f(\eta'))\rho_x(\xi', \eta')}{\rho_x(\xi, \xi')\rho_x(\xi, \eta')
\rho_y(f(\xi'), f(\eta'))}
$$
For a quadruple $(\xi,\xi',\eta,\eta') \in \partial^4 X$, by Theorem \ref{crossratiodistn} we have
$$
e^{-2||S(f)||_{\infty}} \leq \frac{[f(\xi), f(\xi'), f(\eta), f(\eta')]}{[\xi, \xi', \eta, \eta']} \leq e^{2||S(f)||_{\infty}}
$$
Passing to the limit above as $\eta \to \xi$, the term in the middle converges to $df_{\rho_x, \rho_y}(\xi)/\delta(\xi, \xi', \eta')$, thus we
 may write $df_{\rho_x, \rho_y}(\xi) = \delta(\xi, \xi', \eta')\cdot E(\xi, \xi', \eta')$ where
$e^{-2||S(f)||_{\infty}} \leq E(\xi, \xi', \eta') \leq e^{2||S(f)||_{\infty}}$.

\medskip

Now given $(\xi, \eta) \in \partial^2 X$, choose $\beta \in \partial X$ distinct from $\xi, \eta$. Then we have:

\begin{align*}
df_{\rho_x, \rho_y}(\xi)df_{\rho_x,\rho_y}(\eta) & = \delta(\xi, \eta, \beta) \delta(\eta, \xi, \beta) E(\xi, \eta, \beta) E(\eta, \xi, \beta) \\
                                                 & = \left(\frac{\rho_y(f(\xi), f(\eta))}{\rho_x(\xi, \eta)}\right)^2 E(\xi, \eta, \beta) E(\eta, \xi, \beta) \\
\end{align*}

so the Theorem follows since $e^{-4||S(f)||_{\infty}} \leq E(\xi, \eta, \beta) E(\eta, \xi, \beta) \leq e^{4||S(f)||_{\infty}}$. $\diamond$

\medskip

\section{Marked length spectrum, geodesic conjugacies, and
Moebius structure at infinity}

\medskip

The following Theorem follows from results of Bourdon (\cite{bourdon1}),
Hamenstadt (\cite{hamenstadt1}) and Otal (\cite{otal1}), we give a
proof for the benefit of the reader.

\medskip

\begin{theorem} \label{closedneg} (Bourdon, Hamenstadt, Otal) Let $X, Y$ be closed
$n$-dimensional Riemannian manifolds with sectional curvatures bounded
above by $-1$, and let $\tilde{X}, \tilde{Y}$
denote their universal covers. Then the
following are equivalent:

\smallskip

\noindent (1) The marked length spectra of $X$ and $Y$ coincide,
i.e. there is an isomorphism $\Phi : \pi_1(X) \to \pi_1(Y)$ such
that $l_Y \circ \Phi = l_X$.

\smallskip

\noindent (2) There is an equivariant Moebius map $f : \partial
\tilde{X} \to \partial \tilde{Y}$

\smallskip

\noindent (3) There is a homeomorphism $\phi : T^1 X \to T^1 Y$
conjugating the geodesic flows.

\end{theorem}

\medskip

\noindent{\bf Proof:} We prove:

\noindent 1. (1) $\Rightarrow$ (2): It is well known that the
isomorphism $\Phi$ induces an equivariant homeomorphism $f :
\partial \tilde{X} \to \partial \tilde{Y}$ such that $f \circ \gamma = \Phi(\gamma) \circ f$
for $\gamma \in \pi_1(X)$ (with $\pi_1(X), \pi_1(Y)$ identified with groups of homeomorphisms
of $\partial \tilde{X}, \partial \tilde{Y}$).

\medskip

Let $h_X, h_Y$ denote the topological entropies of the geodesic
flows of $X$ and $Y$. For $t \geq 0$, let $\nu_X(t), \nu_Y(t)$ denote the number of
conjugacy classes $[\gamma], [\gamma']$ in $\pi_1(X), \pi_1(Y)$
with $l_X(\gamma) \leq t, l_Y(\gamma') \leq t$. Then by
hypothesis, $\nu_X(t) \equiv \nu_Y(t)$. Hence from Bowen's formula
for the topological entropy (\cite{bowen2}) we have
$$
h_X = \lim_{t \to +\infty} \frac{\log (\nu_X(t))}{t} = \lim_{t \to +\infty} \frac{\log
(\nu_Y(t))}{t} = h_Y
$$
Let $\mu_X, \mu_Y$ denote the Bowen-Margulis currents on
$\partial^2 \tilde{X}, \partial^2 \tilde{Y}$; these are the geodesic currents corresponding to the
Bowen-Margulis measures on $T^1 X, T^1 Y$, the unique invariant measures of maximal entropy. Then it follows from
Bowen's formula for the Bowen-Margulis measure (\cite{bowen2})
that for any fixed small $\epsilon > 0$,
$$
\mu_X = \lim_{t \to +\infty} \frac{1}{N_{\epsilon, X}(t)}
\sum_{[\gamma] \in CO_{\epsilon,X}(t)} \delta_{[\gamma]}
$$
where $CO_{\epsilon,X}(t)$ is the set of conjugacy classes $[\gamma]$ in $\pi_1(X)$
with $l_X(\gamma) \in [t-\epsilon,t+\epsilon]$, $N_{\epsilon,
X}(t)$ is the cardinality of $CO_{\epsilon,X}(t)$, and
$\delta_{[\gamma]}$ denotes the atomic geodesic current associated
to a conjugacy class $[\gamma]$.

\medskip

Since $\Phi$ preserves lengths, it follows that $(f \times f)_*(\mu_X) =
\mu_Y$. We now recall Kaimanovich's formula for the Bowen-Margulis
current (\cite{kaimanovich}),
$$
d\mu_X(\xi,\eta) = \frac{d\nu_{x,X}(\xi)
d\nu_{x,X}(\eta)}{(\rho_x(\xi,\eta))^{2h_X}}
$$
where $x \in \tilde{X}$ (the right-hand side above is independent of the
choice of $x$) and $\nu_{x,X}$ is the Patterson-Sullivan measure
on $\partial \tilde{X}$ based at the point $x$.

\medskip

\noindent{\bf Claim.} For any $x \in \tilde{X}, y \in \tilde{Y}$,
the map $f$ is absolutely continuous with respect to the Patterson-Sullivan measures
$\nu_{x,X}, \nu_{y,Y}$.

\medskip

\noindent{\bf Proof of claim:} Let $A \subset \partial X$ such that $\nu_{x,X}(A) = 0$.
Let $U, V \subset \partial X$ be
closed disjoint balls in $(\partial X, \rho_x)$, let $\delta$ denote the minimum distance
between points of $U$ and $V$. Let $A' = A \cap
U$. Then we have
\begin{align*}
\nu_{y,Y}(f(A')) \nu_{y,Y}(f(V)) & \leq \mu_Y(f(A') \times f(V))
\\
& = \mu_X(A' \times V) \\
& \leq \frac{\nu_{x,X}(A') \nu_{x,X}(V)}{\delta^{2 h_X}} = 0 \\
\end{align*}
hence $\nu_{y,Y}(f(A')) = 0$. It follows that $\nu_{y,y}(f(A)) =
0$. This proves the claim.

\medskip

Let $g$ be the Radon-Nikodym derivative of $f^{-1}_* \nu_{y,Y}$ with
respect to $\nu_{x,X}$. Then the equality $(f \times f)_*(\mu_X) =
\mu_Y$ implies that for $\mu_X$-a.e. $(\xi,\eta) \in \partial^2
\tilde{X}$ we have
$$
\frac{\rho_y(f(\xi), f(\eta))^{2h_Y}}{\rho_x(\xi,
\eta)^{2h_X}} = g(\xi)g(\eta) \ ,
$$
in particular the above equality holds for $(\xi,\eta)$ in a dense
subset $A \subset \partial^2 \tilde{X}$. Since $h_X = h_Y$, it follows that $f$
preserves cross-ratios of quadruples in the dense subset
$\partial^2 A \subset \partial^4 \tilde{X}$, and hence preserves
all cross-ratios, since cross-ratios are continuous.

\medskip

\noindent 2. (2) $\Rightarrow$ (3): Let $\phi : T^1 \tilde{X} \to T^1 \tilde{Y}$
be the geodesic conjugacy induced by $f$, as given by Proposition
\ref{confconj}. Then it is easy to see that $\phi$ is equivariant,
hence induces a geodesic conjugacy $\overline{\phi} : T^1 X \to
T^1 Y$.

\medskip

\noindent 3. (3) $\Rightarrow$ (1): The conjugacy $\phi$ induces
an equivariant conjugacy $\tilde{\phi} : T^1 \tilde{X} \to T^1
\tilde{Y}$, which is uniformly continuous since $\phi$ is
uniformly continuous, hence by Theorem \ref{geodesicmoebius} there
is a Moebius homeomorphism $f : \partial \tilde{X} \to \partial
\tilde{Y}$ such that $\mathcal{E}_Y(\tilde{\phi}(\gamma)) = (f \times f)
\circ \mathcal{E}_X(\gamma)$. Moreover $f$ is equivariant because
$\tilde{\phi}$ is. Identifying $\pi_1(X), \pi_1(Y)$ with groups
of homeomorphisms of $\partial \tilde{X}, \partial \tilde{Y}$, we
obtain a map
\begin{align*}
\Phi : \pi_1(X) & \to \pi_2(Y) \\
          g     & \mapsto f \circ g \circ f^{-1} \\
\end{align*}
which is clearly an isomorphism.

\medskip

Each $g \in \pi_1(X)$ has a unique attracting and
a unique repelling fixed point on $\partial \tilde{X}$, denoted $\xi^+_g, \xi^-_g$ respectively. For
any $\gamma \in \mathcal{G}\tilde{X}$ with $\mathcal{E}_{\tilde{X}}(\gamma) = (\xi^-_g, \xi^+_g)$, we have
$g(\gamma'(0)) = \phi^{\tilde{X}}_{t_1}(\gamma'(0))$, where $t_1 = l_X(g)$. Now
$f(\xi^+_g), f(\xi^-_g)$ are the attracting and repelling fixed
points of $\Phi(g)$, and $\tilde{\phi}(\gamma) \in
\mathcal{G}\tilde{Y}$ satisfies $\mathcal{E}_{\tilde{X}}(\tilde{\phi}(\gamma)) = (f(\xi^-_g),
f(\xi^+_g))$ (we are abusing notation writing $\tilde{\phi}$ also
for the induced map $\mathcal{G}\tilde{X} \to
\mathcal{G}\tilde{Y}$). Hence $\Phi(g)(\tilde{\phi} \circ
\gamma'(0)) = \phi^{\tilde{Y}}_{t_2}(\tilde{\phi} \circ
\gamma'(0))$ where $t_2 = l_Y(\Phi(g))$.

\medskip

Since $\tilde{\phi}$ is equivariant and is a geodesic conjugacy,
we also have
$$
\Phi(g)(\tilde{\phi} \circ \gamma'(0)) = \tilde{\phi}(g(\gamma'(0))) = \tilde{\phi}(\phi^{\tilde{X}}_{t_1}(\gamma'(0))) =
\phi^{\tilde{Y}}_{t_1}(\tilde{\phi} \circ \gamma'(0)).
$$
Since the time-$t$-map of the geodesic flow of $\tilde{Y}$ has no fixed points
for $t \neq 0$, we must have $t_1 = t_2$, i.e. $l_Y(\Phi(g)) = l_X(g)$. $\diamond$

\medskip

We obtain as a corollary the following:

\medskip

\begin{theorem} \label{equivconf} Let $X, Y$ be closed $n$-dimensional
Riemannian manifolds with sectional curvatures bounded above by
$-1$, and let $\tilde{X}, \tilde{Y}$ be their universal covers. If
$f : \partial \tilde{X} \to \partial \tilde{Y}$ is an equivariant
$C^1$ conformal map, then $f$ is Moebius.
\end{theorem}

\medskip

\noindent{\bf Proof:} Let $\phi : T^1 \tilde{X} \to T^1 \tilde{Y}$
be the geodesic conjugacy given by Proposition \ref{confconj}.
Then the equivariance of $f$ implies that of
$\phi$, hence $\phi$ is the lift of a conjugacy $\overline{\phi} :
T^1 X \to T^1 Y$ which is uniformly continuous, hence $\phi$ is
uniformly continuous. It then follows from Theorem
\ref{geodesicmoebius} that $f$ is Moebius. $\diamond$

\medskip

\section{Nearest points and almost isometric extension of Moebius maps}

\medskip

Let $X$ be a proper geodesically complete CAT(-1) space such that
$\partial X$ has at least four points, and let $\mathcal{M} = \mathcal{M}(\partial X, \rho_x)$.
Since the image of the isometric embedding $X \to \mathcal{M}$ is
closed in $\mathcal{M}$ and the space $\mathcal{M}$ is proper, it
follows that for all $\rho \in \mathcal{M}$ there exists $x \in
X$ minimizing $d_{\mathcal{M}}(\rho, \rho_{y})$ over  $y \in X$.

\medskip

\begin{theorem} \label{roottwodense} The image of the map $i_X : X \to \mathcal{M}$ is
$\frac{1}{2}\log 2$-dense in $\mathcal{M}$.
\end{theorem}

\medskip

\noindent{\bf Proof:} Given $\rho \in \mathcal{M}$ let $x \in X$
minimize $d_{\mathcal{M}}(\rho, \rho_y)$ over $y \in X$. Let $\lambda = \sup \log
\frac{d\rho}{d\rho_x} = d_{\mathcal{M}}(\rho, \rho_x)$, let $Z \subset \partial X$ be the set where $\log
\frac{d\rho}{d\rho_x} = \lambda$ and let $\xi_0 \in Z$.

\medskip

Suppose that $\lambda > \frac{1}{2}\log 2$. Then for any $\xi \in Z$, by the Geometric Mean Value Theorem we have
$$
1 \geq \rho(\xi_0,\xi)^2 = \rho_x(\xi_0, \xi)^2 \frac{d\rho}{d\rho_x}(\xi_0) \frac{d\rho}{d\rho_x}(\xi) = \rho_x(\xi_0, \xi)^2 e^{2 \lambda} > \rho_x(\xi_0, \xi)^2 \cdot 2
$$
hence $\max_{\xi \in Z} \rho_x(\xi_0, \xi) < 1/\sqrt{2}$. It follows that there is an open neighbourhood $N \supset Z$ and $\epsilon > 0$ such that
$\angle^{(-1)}\xi x \xi_0 \leq \pi/2 - \epsilon$ for all $\xi \in N$. By monotonicity of comparison angles, for any $y \in [x,\xi_0)$, we also have
$\angle^{(-1)}\xi x y \leq \pi/2 - \epsilon$ for all $\xi \in N$, so $\cos(\angle^{(-1)}\xi x y) \geq \delta_0$ for some $\delta_0 > 0$.
Now let $\lambda' = \sup_{\xi \in \partial X - N} \log
\frac{d\rho}{d\rho_x}(\xi), \delta_1 = \lambda - \lambda' > 0$, then, using Lemma \ref{embedderiv}, let $t_0 < \delta_1/3$ be such that, for
$y \in [x,\xi_0)$ at distance $t$ from $x$, we have
$$
g_{x,y}(\xi) = t \cos(\angle^{(-1)}\xi x y) + o(t)
$$
where $||o(t)||_{\infty} < t \delta_0 / 2$ for $t \leq t_0$. Then, using the Chain Rule,
we have for $\xi \in N$ and $0 < t \leq t_0$, letting $y \in [x, \xi_0)$ be the point at distance $t$ from $x$,
$$
\log \frac{d\rho}{d\rho_{y}}(\xi) = \log \frac{d\rho}{d\rho_x}(\xi) -
g_{x,y}(\xi) \leq \lambda - t \delta_0 + t \delta_0/2 < \lambda,
$$
while for $\xi \in \partial X - N$ and $0 < t \leq t_0$ we have
$$
\log \frac{d\rho}{d\rho_{y}}(\xi) = \log \frac{d\rho}{d\rho_x}(\xi) -
g_{x,y}(\xi) \leq \lambda' + t + t \delta_0 / 2 \leq \lambda - \delta_1 +
2\delta_1/3 < \lambda
$$
hence for $0 < t \leq t_0$ we have $d_{\mathcal{M}}(\rho, \rho_{y}) < d_{\mathcal{M}}(\rho, \rho_x)$, a
contradiction. $\diamond$

\medskip

\begin{theorem} \label{treesurjective} If $X$ is a metric tree
then the map $i_X: X \to \mathcal{M}$ is a surjective isometry.
\end{theorem}

\medskip

\noindent{\bf Proof:} Suppose not, let $\rho \in \mathcal{M}$ be a point not in the
image, let $x \in X$ minimize $d_{\mathcal{M}}(\rho, \rho_y)$ over $y \in X$. Let $\lambda = \sup \log
\frac{d\rho}{d\rho_x} > 0$, let $Z \subset \partial X$ be the set where $\log
\frac{d\rho}{d\rho_x} = \lambda$ and let $\xi_0 \in Z$. Then for all
$\xi \in Z$, we have $1 \geq \rho_x(\xi_0, \xi) e^{\lambda}$ hence
$\rho_x(\xi_0, \xi) \leq e^{-\lambda}$. Let $0 < \lambda' < \lambda$, and choose a
neighbourhood $N \supset Z$ such that $\rho_x(\xi_0, \xi) \leq e^{-\lambda'}$ for all
$\xi \in N$. Letting $y_0$ be the point on
the ray $[x,\xi)$ at distance $\lambda'$ from $a$, since $X$ is a tree it follows that
the segment $[x, y_0]$ is contained in all the rays $[x, \xi), \xi
\in N$. Hence for $0 < t \leq \lambda'$, it follows that $\angle^{(-1)}\xi x y = 0$ for all $\xi \in N$, where
$y$ is the point on $[x, \xi_0)$ at distance $t$ from $x$. Thus $\cos(\angle^{(-1)}\xi x y) = 1$ for all $\xi \in N$, and
now the same argument as in the proof of Theorem \ref{roottwodense} above shows that we may choose $0 < t_0 < \lambda'$ such that
for $0 < t \leq t_0$ we have $d_{\mathcal{M}}(\rho, \rho_{y}) < d_{\mathcal{M}}(\rho, \rho_x)$, a
contradiction. $\diamond$

\medskip

Now let $X, Y$ be proper geodesically complete CAT(-1) spaces such that $\partial X$ has
at least four points, let $f : \partial X \to \partial Y$ be
a Moebius homeomorphism, and let
$\mathcal{M}_X = \mathcal{M}(\partial X, \rho_x), \mathcal{M}_Y = \mathcal{M}(\partial Y,
\rho_y)$ where $x \in X, y \in Y$. Let $g = f^{-1}$, then for $\rho \in \mathcal{M}_X$ we can define
the pull-back metric $g_* \rho$ on $\partial Y$ by $g_* \rho(\xi, \xi') := \rho(g(\xi), g(\xi')),
\xi, \xi' \in \partial Y$. Since $g$ is Moebius it follows easily that $g_* \rho \in \mathcal{M}_Y$.
We can therefore define a map
\begin{align*}
\hat{f} :  \mathcal{M}_X & \to \mathcal{M}_Y \\
              \rho       & \mapsto g_* \rho   \\
\end{align*}
which it is easy to see is a surjective isometry.

\medskip

We define a nearest-point
projection map for $X$,
\begin{align*}
\pi_X : \mathcal{M}_X & \to X \\
             \rho     & \mapsto a \\
\end{align*}
by choosing for each $\rho \in \mathcal{M}_X$ a point
$a \in X$ minimizing $d_{\mathcal{M}_X}(\rho,
\rho_x), x \in X$ (not necessarily unique), and similarly we
define a map $\pi_Y :\mathcal{M}_Y \to Y$.
We can now prove the Theorems \ref{mainthm1} and \ref{mainthm2}:

\medskip

\noindent{\bf Proof of Theorem \ref{mainthm1}:} Define $F: X \to
Y$ by $F = \pi_Y \circ \hat{f} \circ i_X$. Then
by Theorem \ref{roottwodense} for $x, x' \in X$, letting $y = F(x), y' = F(x')$ we have
\begin{align*}
|d_Y(y, y') - d_X(x, x')| & = |d_{\mathcal{M}_Y}(\rho_y, \rho_y') -
d_{\mathcal{M}_Y}(\hat{f}(\rho_x), \hat{f}(\rho_{x'}))| \\
& \leq |d_{\mathcal{M}_Y}(\rho_y, \rho_{y'}) -
d_{\mathcal{M}_Y}(\hat{f}(\rho_x), \rho_{y'})| + |d_{\mathcal{M}_Y}(\hat{f}(\rho_x), \rho_{y'}) -
d_{\mathcal{M}_Y}(\hat{f}(\rho_x), \hat{f}(\rho_{x'}))| \\
& \leq d_{\mathcal{M}_Y}(\rho_y, \hat{f}(\rho_x)) + d_{\mathcal{M}_Y}(\rho_{y'},
\hat{f}(\rho_{x'})) \leq \log 2 \\
\end{align*}
so $F$ is a $(1, \log 2)$-quasi-isometry. Given $y \in Y$, by Theorem \ref{roottwodense} we may choose
$x \in X$ such that $d_{\mathcal{M}_X}(f_* \rho_y, \rho_x) \leq
\frac{1}{2} \log 2$, then by definition of $F$,

\begin{align*}
d_Y(F(x), y) = d_{\mathcal{M}_Y}(\rho_{F(x)}, \rho_y) & \leq d_{\mathcal{M}_Y}(\rho_{F(x)}, \hat{f}(\rho_x)) + d_{\mathcal{M}_Y}(\hat{f}(\rho_x), \rho_y) \\
                                             & \leq d_{\mathcal{M}_Y}(\rho_y, \hat{f}(\rho_x)) + d_{\mathcal{M}_Y}(\hat{f}(\rho_x), \rho_y) \\
                                             & = 2 d_{\mathcal{M}_X}(f_* \rho_y, \rho_x) \\
                                             & \leq \log 2 \\
\end{align*}

thus the image of $F$ is $\log 2$-dense in $Y$. 

\medskip

It follows from the above that $F$ has a continuous extension $\partial F : \partial X \to \partial Y$,
it remains to prove that $\partial F = f$. Let $\xi \in \partial
X, x \in X$ and let $a \in X$ converge to $\xi$ along the ray $[x,
\xi)$. Let $y = F(x), b = F(a), \lambda = d_Y(y,b)$, then $b \to \eta = \partial
F(\xi), \lambda \geq d_X(x,a) - \log 2 \to \infty$ as $a \to \xi$.
Extend $[y, b]$ to a geodesic ray $[y,\eta')$ where $\eta' \in \partial Y$,
then $\frac{d\rho_b}{d\rho_y}(\eta') = e^{\lambda}$ and $b \to \eta$ implies $\eta' \to \eta$. By the Chain Rule,
$$
||\log \frac{d \rho_b}{d \rho_y} - \log \frac{d g_* \rho_a}{d g_*
\rho_x}||_{\infty} \leq d_{\mathcal{M}_Y}(\rho_b, g_* \rho_a) + d_{\mathcal{M}_Y}(\rho_y, g_* \rho_x) \leq \log 2
$$
and $\log \frac{d g_* \rho_a}{d g_* \rho_x}(f(\xi)) = d_X(x,a) \geq
d_Y(y,b) - \log 2$, hence
$$
\log \frac{d \rho_b}{d \rho_y}(f(\xi)) \geq \log \frac{d g_* \rho_a}{d g_* \rho_x}(f(\xi)) - \log 2 \geq \lambda
- 2 \log 2
$$
so $\frac{d \rho_b}{d \rho_y}(f(\xi)) \geq e^{\lambda}/4$, thus
$$
1 \geq \rho_b(f(\xi), \eta')^2 = \rho_y(f(\xi), \eta')^2 \frac{d \rho_b}{d \rho_y}(f(\xi)) \frac{d \rho_b}{d \rho_y}(\eta') \geq
\rho_y(f(\xi), \eta')^2 e^{2\lambda}/4
$$
hence $\rho_y(f(\xi), \eta') \to 0$, and $\eta' \to \eta$, so $f(\xi) = \eta = \partial
F(\xi)$. $\diamond$

\medskip

\noindent{\bf Proof of Theorem \ref{mainthm2}:} For $X, Y$ proper
geodesically complete metric trees such that $\partial X$ has at least four points, by Theorem
\ref{treesurjective} we have surjective isometries $i_X : X \to
\mathcal{M}_X, \hat{f} : \mathcal{M}_X \to \mathcal{M}_Y, i^{-1}_Y : \mathcal{M}_Y \to Y$, and
it is clear that the map $F$ defined above equals the composition of these isometries,
hence is a surjective isometry $X \to Y$ extending $f$. $\diamond$

\medskip

\noindent{\bf Proof of Theorem \ref{geodconj}:} The assertion (1)
follows immediately from Theorem \ref{geodesicmoebius} and
\ref{mainthm1}. For the assertion (2), Theorem
\ref{geodesicmoebius} and Theorem \ref{bourdonthm} give us an
isometry $F : X \to Y$ with $f = \partial F$ a Moebius homeomorphism.
Given $y \in Y$, choose a bi-infinite geodesic $\gamma' \in
\mathcal{G}Y$ with $y \in \gamma'(\mathbb{R})$, let $\gamma \in
\mathcal{G}X$ be a geodesic whose endpoints map to those of
$\gamma'$ under $f$, then $F$ maps the image of $\gamma$ onto the
image of $\gamma'$, in particular $y$ belongs to the image of $F$,
hence $F$ is surjective. $\diamond$

\medskip

Finally we prove Theorem \ref{confextn} on almost isometric extension of $C^1$ conformal maps with bounded
integrated Schwarzian. The proof proceeds along similar lines to the proof of Theorem \ref{mainthm1}.

\medskip

Let $(Z, \rho_0)$ be a compact metric space. We assume $Z$ has no isolated points, and that $\rho_0$ is diameter one and antipodal. We define the set of metrics
$$
Conf(Z, \rho_0) := \{ \rho | \rho \hbox{ is a diameter one antipodal metric on } Z \hbox{ s.t. } id: (Z, \rho_0) \to (Z, \rho) \hbox{ is } C^1 \hbox{ conformal} \}
$$

\medskip

Note that $\mathcal{M}(Z, \rho_0) \subset Conf(Z, \rho_0)$. For $\rho_1,\rho_2 \in Conf(Z, \rho_0)$, the derivative $\frac{d\rho_2}{d\rho_1}$ is a
continuous function on $Z$ so we can define
$$
d_{Conf}(\rho_1,\rho_2) := \max_{\xi \in Z} \left|\log \frac{d\rho_2}{d\rho_1}(\xi)\right|
$$

\medskip

Then it is easy to see that $d_{Conf}$ is a pseudo-metric on $Conf(Z, \rho_0)$ (though not necessarily a metric) extending the metric
$d_{\mathcal M}$ on $\mathcal{M}(Z, \rho_0)$. Any $C^1$ conformal map between compact metric spaces $f : Z_1 \to Z_2$ induces a natural
bijective isometry of pseudo-metric spaces $\hat{f} : Conf(Z_1) \to Conf(Z_2)$ by push-forward of metrics.

\medskip

Now let $X$ be a simply connected complete Riemannian manifold with sectional curvatures satisfying $-b^2 \leq K \leq -1$,
let $Y$ be a proper geodesically complete CAT(-1) space and let $f : \partial X \to \partial Y$ be a $C^1$ conformal map with
bounded integrated Schwarzian. We let
$Conf(\partial X) = Conf(\partial X, \rho_x), Conf(\partial Y) = Conf(\partial Y, \rho_y)$ for some $x \in X, y \in Y$ (note
the definition does not depend on the choice of $x$ and $y$), and let $\hat{f} : Conf(\partial X) \to Conf(\partial Y)$ be the induced
isometry. We note that
$$
\frac{d\hat{f}(\rho_x)}{d\rho_y} \circ f = 1/df_{\rho_x, \rho_y}
$$
for all $x \in X, y \in Y$.

\begin{lemma}\label{confmaxmin} For all $x \in X, y \in Y$,
$$
\left|\min_{\xi \in \partial X} \log df_{\rho_x,\rho_y}(\xi) + \max_{\xi \in \partial X} \log df_{\rho_x, \rho_y}(\xi)\right| \leq 4||S(f)||_{\infty}
$$
Moreover
$$
\max_{\xi \in \partial X} |\log df_{\rho_x, \rho_y}(\xi)| \leq - \min_{\xi \in \partial X} \log df_{\rho_x, \rho_y}(\xi) + 4||S(f)||_{\infty}
$$
\end{lemma}

\medskip

\noindent{\bf Proof:} Let $\lambda = \max_{\xi \in \partial X} \log df_{\rho_x, \rho_y}(\xi), \mu = \min_{\xi \in \partial X} \log df_{\rho_x, \rho_y}(\xi)$. Let $\eta \in \partial X$ minimize $\log df_{\rho_x, \rho_y}$. Choose $\xi \in \partial X$ such that $\rho_y(f(\xi), f(\eta)) = 1$,
then we have, using Theorem \ref{confgmvt},

\begin{align*}
e^{\lambda} e^{\mu} & \geq df_{\rho_x, \rho_y}(\xi) df_{\rho_x, \rho_y}(\eta) \\
                    & \geq \left(\frac{\rho_y(f(\xi), f(\eta))}{\rho_x(\xi, \eta)}\right)^2 e^{-4||S(f)||_{\infty}} \\
                    & \geq e^{-4||S(f)||_{\infty}} \\
\end{align*}

so $\lambda + \mu \geq -4||S(f)||_{\infty}$. For the other inequality, let $\eta \in \partial X$ maximize $\log df_{\rho_x, \rho_y}$, choose
$\xi \in \partial X$ such that $\rho_x(\xi, \eta) = 1$, then again by Theorem \ref{confgmvt}, we have

\begin{align*}
e^{\lambda} e^{\mu} & \leq df_{\rho_x, \rho_y}(\xi) df_{\rho_x, \rho_y}(\eta) \\
                    & \leq \left(\frac{\rho_y(f(\xi), f(\eta))}{\rho_x(\xi, \eta)}\right)^2 e^{4||S(f)||_{\infty}} \\
                    & \leq e^{4||S(f)||_{\infty}} \\
\end{align*}

This proves the first assertion above. For the second, let $L = \max_{\xi \in \partial X} |\log df_{\rho_x, \rho_y}(\xi)|$. Then
either $L = -\mu$ or $L = \lambda \leq -\mu + 4||S(f)||_{\infty}$ by the first assertion. $\diamond$

\medskip

\begin{lemma}\label{confdense} For all $x \in X$, there exists $y \in Y$ such that $d_{Conf}(\hat{f}(\rho_x), \rho_y) \leq \frac{1}{2}\log 2
+ 6||S(f)||_{\infty}$.
\end{lemma}

\medskip

\noindent{\bf Proof:} Given $x \in X$, define the function $\phi : Y \to \mathbb{R}$ by $\phi(y) = \max_{\xi \in \partial Y} \log \frac{d\hat{f}(\rho_x)}{d\rho_y}(\xi)$. Note $\phi$ is $1$-Lipschitz (since $i_Y : Y \to Conf(\partial Y)$ is an isometry).
Let $y_n \in Y$ be a sequence such that $\phi(y_n) \to \inf_{y \in Y} \phi(y)$. Then by Lemma \ref{confmaxmin},

\begin{align*}
d_{Conf}(\rho_{y_n}, \hat{f}(\rho_x)) & = \max_{\xi \in \partial X} |\log df_{\rho_x,\rho_{y_n}}(\xi)| \\
                                      & \leq -\min_{\xi \in \partial X} \log df_{\rho_x,\rho_{y_n}}(\xi) + 4||S(f)||_{\infty} \\
                                      & = \max_{\xi \in \partial X} (-\log df_{\rho_x,\rho_{y_n}}(\xi)) + 4||S(f)||_{\infty} \\
                                      & = \max_{\xi \in \partial Y} \log \frac{d\hat{f}(\rho_x)}{d\rho_{y_n}}(\xi) + 4||S(f)||_{\infty} \\
                                      & = \phi(y_n) + 4||S(f)||_{\infty} \\
\end{align*}

Since the sequence $\{\phi(y_n)\}$ is bounded above, by the triangle inequality $d_{Conf}(\rho_{y_n}, \rho_{y_m})$ is bounded independent of $m,n$, hence
so is $d_Y(y_n, y_m)$. Thus we have a convergent subsequence $y_{n_k} \to z \in Y$, and $\phi(z) = \lim \phi(y_{n_k}) = \inf_{y \in Y} \phi(y)$.

\medskip

\noindent{\bf Claim.} Let $\lambda = \phi(z)$, then $\lambda \leq \frac{1}{2}\log 2 + 2||S(f)||_{\infty}$.

\medskip

\noindent{\bf Proof of Claim:} Suppose $\lambda > \frac{1}{2}\log 2 + 2||S(f)||_{\infty}$. Let $Z \subset \partial Y$ be the set where
$\log \frac{d\hat{f}(\rho_x)}{d\rho_{z}} = \lambda$, and let $\xi_0 \in Z$. Then for any $\xi \in Z$, by Theorem \ref{confgmvt} we have:

\begin{align*}
1 & \geq \hat{f}(\rho_x)(\xi_0, \xi)^2 \\
  & \geq \rho_{z}(\xi_0, \xi)^2 \frac{d\hat{f}(\rho_x)}{d\rho_{z}}(\xi_0) \frac{d\hat{f}(\rho_x)}{d\rho_{z}}(\xi) e^{-4||S(f)||_{\infty}} \\
  & = \rho_{z}(\xi_0, \xi)^2 e^{2\lambda} e^{-4||S(f)||_{\infty}} \\
  & > 2 \rho_z(\xi_0, \xi)^2
\end{align*}

thus $\rho_z(\xi_0, \xi) < 1/\sqrt{2}$. It follows that there is an open neighbourhood $N \supset Z$ and $\epsilon > 0$ such that
$\angle^{(-1)}\xi z \xi_0 \leq \pi/2 - \epsilon$ for all $\xi \in N$. By monotonicity of comparison angles, for any $y \in [z,\xi_0)$, we also have
$\angle^{(-1)}\xi z y \leq \pi/2 - \epsilon$ for all $\xi \in N$, so $\cos(\angle^{(-1)}\xi z y) \geq \delta_0$ for some $\delta_0 > 0$.
Now let $\lambda' = \sup_{\xi \in \partial Y - N} \log
\frac{d\hat{f}(\rho_x)}{d\rho_z}(\xi), \delta_1 = \lambda - \lambda' > 0$, then, using Lemma \ref{embedderiv}, let $t_0 < \delta_1/3$ be such that, for
$y \in [z,\xi_0)$ at distance $t$ from $z$, we have
$$
g_{z,y}(\xi) = t \cos(\angle^{(-1)}\xi z y) + o(t)
$$
where $||o(t)||_{\infty} < t \delta_0 / 2$ for $t \leq t_0$. Then, using the Chain Rule,
we have for $\xi \in N$ and $0 < t \leq t_0$, letting $y \in [z, \xi_0)$ be the point at distance $t$ from $z$,
$$
\log \frac{d\hat{f}(\rho_x)}{d\rho_y}(\xi) = \log \frac{d\hat{f}(\rho_x)}{d\rho_z}(\xi) -
g_{z,y}(\xi) \leq \lambda - t \delta_0 + t \delta_0/2 < \lambda,
$$
while for $\xi \in \partial X - N$ and $0 < t \leq t_0$ we have
$$
\log \frac{d\hat{f}(\rho_x)}{d\rho_y}(\xi) = \log \frac{d\hat{f}(\rho_x)}{d\rho_z}(\xi) -
g_{z,y}(\xi) \leq \lambda' + t + t \delta_0 / 2 \leq \lambda - \delta_1 +
2\delta_1/3 < \lambda
$$
hence for $0 < t \leq t_0$ we have $\phi(y) < \phi(z)$, a
contradiction. This proves the Claim.

\medskip

Now it follows from Lemma \ref{confmaxmin} that
\begin{align*}
d_{Conf}(\hat{f}(\rho_x), \rho_z) & = \max_{\xi \in \partial Y} |\log df_{\rho_x, \rho_z}(\xi)| \\
                                  & \leq -\min_{\xi \in \partial Y} \log df_{\rho_x, \rho_z}(\xi) + 4||S(f)||_{\infty} \\
                                  & = \lambda + 4||S(f)||_{\infty} \\
                                  & \leq \frac{1}{2}\log 2 + 6||S(f)||_{\infty} \\
\end{align*}

$\diamond$

\medskip

We can now prove Theorem \ref{confextn}:

\medskip

\noindent{\bf Proof:} By the same argument as in the previous Lemma, for each $x \in X$ we may choose a point $F(x) \in Y$ which minimizes
$d_{Conf}(\hat{f}(\rho_x), \rho_{y})$ over $y \in Y$, and
moreover we have $d_{Conf}(\hat{f}(\rho_x), \rho_{F(x)}) \leq \frac{1}{2}\log 2 + 6||S(f)||_{\infty}$. This defines a map $F : X \to Y$.

\medskip

For $p,q \in X$, let $u = F(p), v = F(q)$, then we have

\begin{align*}
|d_Y(u, v) - d_X(p,q)| & = |d_{Conf}(\rho_u, \rho_v) - d_{Conf}(\hat{f}(\rho_p), \hat{f}(\rho_q))| \\
                       & \leq d_{Conf}(\rho_u, \hat{f}(\rho_p)) + d_{Conf}(\rho_v, \hat{f}(\rho_q)) \\
                       & \leq \log 2 + 12||S(f)||_{\infty} \\
\end{align*}

thus $F$ is a $(1, \log 2 + 12||S(f)||_{\infty})$-quasi-isometry.

\medskip

Thus $F$ has a continuous extension to the boundary $\partial F : \partial X \to \partial Y$.

\medskip

We prove $\partial F = f$. Let $\xi \in \partial
X, x \in X$ and let $a \in X$ converge to $\xi$ along the ray $[x,
\xi)$. Let $y = F(x), b = F(a), \lambda = d_Y(y,b)$, then $b \to \eta = \partial
F(\xi), \lambda \geq d_X(x,a) - \log 2 - 12||S(f)||_{\infty} \to \infty$ as $a \to \xi$.
Extend $[y, b]$ to a geodesic ray $[y,\eta')$ where $\eta' \in \partial Y$,
then $\frac{d\rho_b}{d\rho_y}(\eta') = e^{\lambda}$ and $b \to \eta$ implies $\eta' \to \eta$. By the Chain Rule,
$$
\left|\log \frac{d \rho_b}{d \rho_y} - \log \frac{d \hat{f}(\rho_a)}{d \hat{f}(
\rho_x)}\right| \leq d_{Conf}(\rho_b, \hat{f}(\rho_a)) + d_{Conf}(\rho_y, \hat{f}(\rho_x)) \leq \log 2 + 12||S(f)||_{\infty}
$$
and $\log \frac{d \hat{f}(\rho_a)}{d \hat{f}(\rho_x)}(f(\xi)) = d_X(x,a) \geq
d_Y(y,b) - \log 2 - 12||S(f)||_{\infty}$, hence
$$
\log \frac{d \rho_b}{d \rho_y}(f(\xi)) \geq \log \frac{d \hat{f}(\rho_a)}{d \hat{f}(\rho_x)}(f(\xi)) - \log 2 - 12||S(f)||_{\infty} \geq \lambda
- 2 \log 2 - 24||S(f)||_{\infty}
$$
so $\frac{d \rho_b}{d \rho_y}(f(\xi)) \geq Ce^{\lambda}$ for some constant $C > 0$, thus
$$
1 \geq \rho_b(f(\xi), \eta')^2 = \rho_y(f(\xi), \eta')^2 \frac{d \rho_b}{d \rho_y}(f(\xi)) \frac{d \rho_b}{d \rho_y}(\eta') \geq
\rho_y(f(\xi), \eta')^2 Ce^{2\lambda}
$$
hence $\rho_y(f(\xi), \eta') \to 0$, so $f(\xi) = \eta = \partial
F(\xi)$.

\medskip

Finally, if $Y$ is also a simply connected complete Riemannian manifold with sectional curvatures satisfying $-b^2 \leq K \leq -1$,
then, given $y \in Y$, we may apply Lemma \ref{confdense} to the map $f^{-1}$ to obtain $x \in X$ such that $d_{Conf}(\hat{f}^{-1}(\rho_y), \rho_x)
\leq \frac{1}{2}\log + 6||S(f)||_{\infty}$ (note $||S(f^{-1})||_{\infty} = ||S(f)||_{\infty}$). Then by definition of $F$,

\begin{align*}
d_Y(F(x), y) = d_{Conf}(\rho_{F(x)}, \rho_y) & \leq d_{Conf}(\rho_{F(x)}, \hat{f}(\rho_x)) + d_{Conf}(\hat{f}(\rho_x), \rho_y) \\
                                             & \leq d_{Conf}(\rho_y, \hat{f}(\rho_x)) + d_{Conf}(\hat{f}(\rho_x), \rho_y) \\
                                             & = 2 d_{Conf}(\hat{f}^{-1}(\rho_y), \rho_x) \\
                                             & \leq \log 2 + 12||S(f)||_{\infty} \\
\end{align*}

thus the image of $F$ is $\log 2 + 12||S(f)||_{\infty}$-dense in $Y$. $\diamond$

\medskip

\section{Dynamical classification of Moebius self-maps}

\medskip

Let $X$ be a proper geodesically complete CAT(-1) space whose
boundary has at least four points. We use the results of the
previous section to prove the dynamical classification of Moebius
self-maps of $\partial X$ stated in Theorem \ref{classfn}:

\medskip

\noindent {\bf Proof of Theorem \ref{classfn}}: Let $f : \partial
X \to \partial X$ be a Moebius homeomorphism. As in the previous section choose
and fix a nearest point projection $\pi_X :
\mathcal{M}(\partial X) \to X$, so for all $\rho \in
\mathcal{M}(\partial X)$, the visual metric $\rho_{x_0}$, where $x_0 = \pi(\rho)$,
minimizes $d_{\mathcal{M}}(\rho, \rho_x), x \in X$. Note in
particular that $\pi_X$ is a $(1, \log 2)$-quasi-isometry, $\pi_X \circ i_X = id_X$
and $d_{\mathcal{M}}(\rho, i_X \circ \pi_X (\rho)) \leq \frac{1}{2}
\log 2, \rho \in \mathcal{M}(\partial X)$, i.e. $i_X \circ \pi_X$ is at a
uniformly bounded distance from $id_{\mathcal{M}(\partial X)}$.

\medskip

Define as in the proof of Theorem \ref{mainthm1} a sequence of
$(1, \log 2)$-quasi-isometric extensions $(F_n : X \to X)_{n \in \mathbb{Z}}$ of the maps $(f^n :
\partial X \to \partial X)_{n \in \mathbb{Z}}$ by putting $F_n = \pi_X \circ \widehat{f^n}
\circ i_X$ where $\widehat{f^n} : \mathcal{M}(\partial X) \to
\mathcal{M}(\partial X)$ denotes the isometry induced by $f^n$.
Note that $\widehat{f^n} = \hat{f}^n$ and $F_0 = id_X$. It is easy to see that
since $i_X \circ \pi_X$ is at a bounded distance from
$id_{\mathcal{M}(\partial X)}$, for any $m,n \in \mathbb{Z}$ the
maps $F_m \circ F_n = \pi_X \circ \widehat{f^m} \circ (i_X \circ \pi_X)
\circ \widehat{f^n} \circ i_X, F_n \circ F_m = \pi_X \circ \widehat{f^n} \circ (i_X \circ \pi_X)
\circ \widehat{f^m} \circ i_X$ and $F_{m+n} = \pi_X \circ \widehat{f^{m+n}}
\circ i_X$ are all within bounded distance of each other.

\medskip

We note that by the definition of $F_n$, for any $x \in X$, the
maps $f^n : (\partial X, \rho_x) \to (\partial X, \rho_{F_n(x)})$
are uniformly $\sqrt{2}$-bi-Lipschitz.

\medskip

Since the maps $F_n$ are uniform $(1, \log 2)$-quasi-isometries,
it is clear that the set of accumulation points in $\partial X$ of
a sequence $(F_n(x))_{n \in \mathbb{Z}}$ is independent of the
choice of $x \in X$. We denote this set by $\Lambda$. We observe that
if $\xi \in \Lambda$, then there is a sequence $(n_k)$ such that for any
$x \in X$, $F_{n_k}(x) \to \xi$, in particular $F_{n_k}(F_1(x)) \to \xi$,
hence $F_1(F_{n_k}(x)) \to \xi$ (as the two sequences are within bounded distance of each other),
and since $F_1$ has boundary value $f$, it follows that $F_1(F_{n_k}(x)) \to f(\xi)$, hence
$\xi = f(\xi)$. Thus all points of $\Lambda$ are fixed points of $f$. We now consider
three cases:

\medskip

\noindent Case 1. $\Lambda = \emptyset$: Then for any $x \in X$,
the sequence $(F_n(x))_{n \in \mathbb{Z}}$ is bounded, so the
metrics $\rho_x$ and $\rho_{F_n(x)}$ are uniformly bi-Lipschitz to
each other independent of $n$, and it follows from the observation made above that
the maps $f^n : (\partial X, \rho_x) \to (\partial X, \rho_x)$ are
uniformly bi-Lipschitz, so we are in Case 1 of Theorem
\ref{classfn}, the elliptic case.

\medskip

\noindent Case 2. $\Lambda = \{\xi_0\}$: Then for any $x \in X$,
$F_n(x) \to \xi_0$ as $|n| \to +\infty$. We claim that $f^n(\xi)
\to \xi_0$ as $|n| \to +\infty$ for all $\xi \in \partial X$, i.e.
we are in Case 2 of Theorem \ref{classfn}, the parabolic case.

\medskip

Suppose not, then there is a $\xi \neq \xi_0$ such that some subsequence
$f^{n_k}(\xi)$ converges to a $\xi_1 \neq \xi_0$. Fix $x \in X$ belonging to the geodesic
$\gamma = (\xi_0, \xi)$. The images $F_{n_k}(\gamma)$ are uniform $(1,
\log 2)$-quasi-geodesics with endpoints $\xi_0, f^{n_k}(\xi)$,
with the endpoints $f^{n_k}(\xi)$ uniformly bounded away from
$\xi_0$, hence there is a ball $B$ of fixed radius around $x$ such that
$F_{n_k}(\gamma)$ intersects $B$ for all $k$. Choose for each $k$
a point $y_k \in F_{n_k}(\gamma) \cap B$. Then $d(y_k, F_{n_k}(x))
\to +\infty$ as $k \to +\infty$. The distances $d(y_k,
F_{n_k}(x)), d(F_{-n_k}(y_k), F_{-n_k}(F_{n_k}(x)))$ differ by a
uniformly bounded amount (since $F_{-n_k}$'s are uniform
quasi-isometries), as do the distances $d(F_{-n_k}(y_k), F_{-n_k}(F_{n_k}(x)))$,
$d(F_{-n_k}(y_k), x)$ (since the maps $F_{-n} \circ F_n$ are within uniformly
bounded distance of the identity), hence $d(F_{-n_k}(y_k), x) \to
+\infty$.

\medskip

The horospherical distances $B(\xi_0, F_{-n_k}(y_k), x)$,
$B(\xi_0, F_{n_k}(F_{-n_k}(y_k)), F_{n_k}(x))$ differ by a
uniformly bounded amount (since the maps $F_{n_k}$
are uniform quasi-isometries with boundary maps $f^{n_k}$ fixing
$\xi_0$), as do $B(\xi_0, F_{n_k}(F_{-n_k}(y_k))$, $F_{n_k}(x)), B(\xi_0, y_k, F_{n_k}(x))$
(since the maps $F_{-n} \circ F_n$ are within uniformly
bounded distance of the identity), and clearly $B(\xi_0, y_k, F_{n_k}(x)) \to +\infty$,
hence $B(\xi_0, F_{-n_k}(y_k), x) \to
+\infty$. Since the points $F_{-n_k}(y_k)$ lie on uniform quasi-geodesics
$F_{-n_k} \circ F_{n_k}(\gamma)$ with fixed endpoints $\xi_0,
\xi$ and $d(F_{-n_k}(y), x) \to +\infty$, it follows that $F_{-n_k}(y_k) \to
\xi$. Since the points $y_k$ are within uniformly bounded distance of $x$
and the maps $F_{-n_k}$ are uniform quasi-isometries, it follows
that $F_{-n_k}(x) \to \xi$, a contradiction.

\medskip

\noindent Case 3. The set $\Lambda$ has at least two points: Then
pick two distinct points $\xi_+, \xi_- \in \Lambda$, and fix a point $x$ on the
geodesic $\gamma = (\xi_+, \xi_-)$. We may assume
(replacing $f$ by $f^{-1}$ if necessary) that
there is a subsequence $F_{n_k}(x) \to \xi_+$ with $n_k \to
+\infty$.

\medskip

We claim first that $f^{n_k}(\xi) \to \xi_+$ for all $\xi \in
\partial X - \{\xi_-\}$. If not, then there is a $\xi \neq \xi_-$
such that, after passing to a further subsequence if necessary,
the distances $\rho_x(\xi_+, f^{n_k}(\xi))$ are bounded below by a
constant $\epsilon > 0$. Since the points $F_{n_k}(x)$ converge to $\xi_+$
and lie on uniform quasi-geodesics $F_{n_k}(\gamma)$ with fixed endpoints
$\xi_+, \xi_-$, it follows that $\rho_{F_{n_k}(x)}(\xi_-, f^{n_k}(\xi)) \to 0$.
However the maps $f^{n_k} : (\partial X, \rho_x) \to (\partial X,
\rho_{F_{n_k}(x)})$ are uniformly bi-Lipschitz, hence the sequence $\rho_{F_{n_k}(x)}(\xi_-,
f^{n_k}(\xi))$ is bounded below by a positive constant times
$\rho_x(\xi_-, \xi)$, and does not tend to zero, a contradiction.

\medskip

We now claim that $f^{n}(\xi) \to \xi_+$ for all $\xi \in
\partial X - \{\xi_-\}$ as $n \to +\infty$. Denoting by
$df_{p,q}(\xi)$ the derivative of the conformal map $f : (\partial X, \rho_p) \to (\partial X,
\rho_q)$ at a point $\xi \in \partial X$, we have

\begin{align*}
(df_{x,x}(\xi^+))^{n_k} & = df^{n_k}_{x,x}(\xi_+) \\
                        & = df^{n_k}_{x, F_{n_k}(x)}(\xi_+) \cdot
                        \frac{d\rho_x}{d\rho_{F_{n_k}(x)}}(\xi_+)
                        \\
                        & \to 0 \\
\end{align*}

since $df^{n_k}_{x, F_{n_k}(x)}(\xi_+)$ is bounded above (by
$\sqrt{2}$) and $\frac{d\rho_x}{d\rho_{F_{n_k}(x)}}(\xi_+) \to 0$
(as the points $F_{n_k}(x)$ converge to $\xi_+$ and lie along
uniform quasi-geodesics $F_{n_k}(\gamma)$ with fixed endpoints
$\xi_+, \xi_-$). It follows that $df_{x,x}(\xi_+) < 1$, hence
there is a neighbourhood $U$ of $\xi_+$ such that $f^n(\xi) \to
\xi_+$ as $n \to +\infty$ for all $\xi \in U$. Now given $\xi \in
\partial X - \{\xi_-\}$, there is a $k$ such that $f^{n_k}(\xi)
\in U$, hence it follows that $f^n(\xi) \to \xi_+$ as $n \to
+\infty$.

\medskip

Now there is a sequence of integers $m_k$ with $|m_k| \to +\infty$
such that $F_{m_k}(x) \to \xi_-$. By the argument given above, we must have
$m_k \to -\infty$ (otherwise there would be a sequence of positive
integers tending to infinity with $f^n$ converging pointwise on
$\partial X - \{\xi_+\}$ to $\xi_-$, contradicting the conclusion
of the previous paragraph). It follows from the same argument as
above that $f^n(\xi) \to \xi_-$ as $n \to -\infty$ for all $\xi \in \partial X -
\{\xi_+\}$. Hence we are in Case 3 of Theorem \ref{classfn}, the
hyperbolic case. $\diamond$

\bibliography{moeb}
\bibliographystyle{alpha}

\end{document}